\documentclass[11pt]{article}

\usepackage{amsmath,amsfonts,amsthm,amssymb,graphicx,tikz,tikz-cd,url,hyperref}
\usepackage[all,2cell,ps]{xy}

\theoremstyle{plain}
\newtheorem{thm}{Theorem}[section]

\newtheorem{lem}[thm]{Lemma}
\newtheorem{prop}[thm]{Proposition}

\newtheorem{cor}[thm]{Corollary}

\newtheorem{qtn}[thm]{Question}

\theoremstyle{definition}

\newtheorem{rem}[thm]{Remark}

\theoremstyle{remark}

\newcommand{\bbB}{\mathbb{B}}
\newcommand{\bbC}{\mathbb{C}}

\newcommand{\bbH}{\mathbb{H}}

\newcommand{\bbP}{\mathbb{P}}
\newcommand{\bbQ}{\mathbb{Q}}

\newcommand{\bbZ}{\mathbb{Z}}

\newcommand{\calC}{\mathcal{C}}
\newcommand{\calD}{\mathcal{D}}
\newcommand{\calE}{\mathcal{E}}

\newcommand{\calL}{\mathcal{L}}

\newcommand{\al}{\alpha}

\newcommand{\Gam}{\Gamma}

\newcommand{\Del}{\Delta}

\newcommand{\lam}{\lambda}
\newcommand{\Lam}{\Lambda}

\newcommand{\Sig}{\Sigma}

\DeclareMathOperator{\GL}{GL}

\DeclareMathOperator{\PU}{PU}

\DeclareMathOperator{\Id}{Id}

\DeclareMathOperator{\Aut}{Aut}

\DeclareMathOperator{\Stab}{Stab}

\newcommand{\bs}{\backslash}

\newcommand{\lra}{\longrightarrow}

\newcommand{\ssm}{\smallsetminus}

\newcommand{\wh}{\widehat}

\newenvironment{pf}{\begin{proof}}{\end{proof}}

\allowdisplaybreaks

\makeatletter
\let\@@pmod\pmod
\DeclareRobustCommand{\pmod}{\@ifstar\@pmods\@@pmod}
\def\@pmods#1{\mkern4mu({\operator@font mod}\mkern 6mu#1)}
\makeatother

\usetikzlibrary{arrows}

\title{Products of curves as ball quotients}
\author{Matthew Stover\footnote{This material is based upon work supported by Grant Number DMS-2203555 from the National Science Foundation.} \\ \small{Temple University}\\ \small{\textsf{mstover@temple.edu}}}
\date{\today}

\begin{document}

\maketitle

\begin{abstract}
For any $g_1, g_2 \ge 0$, this paper shows that there is a cocompact lattice $\Gam < \PU(2,1)$ such that the ball quotient $\Gam \bs \bbB^2$ is birational to a product $C_1 \times C_2$ of smooth projective curves $C_j$ of genus $g_j$. The only prior examples were $\bbP^1 \times \bbP^1$, due to Deligne--Mostow and rediscovered by many others, and a lesser-known product of elliptic curves whose existence follows from work of Hirzebruch. Combined with related new examples, this answers the rational variant of a question of Gromov in the positive for surfaces of Kodaira dimension $\kappa \le 0$, namely that they admit deformations $V^\prime$ such that there is a compact ball quotient $\Gam \bs \bbB^2$ with a rational map $\Gam \bs \bbB^2 \dashrightarrow V^\prime$. Often the proof gives the stronger conclusion that $V^\prime$ is birational to a ball quotient orbifold. It also follows that every simply connected $4$-manifold is dominated by a complex hyperbolic manifold. All examples considered in this paper are shown to be arithmetic, and even arithmeticity of Hirzebruch's example appears to be new.
\end{abstract}

\section{Introduction}\label{sec:Intro}

If $\Gam < \PU(2,1)$ is a discrete subgroup acting cocompactly on the unit ball $\bbB^2$ in $\bbC^2$, then the ball quotient $\Gam \bs \bbB^2$ is isomorphic to a normal projective variety $X$. A realization of a normal projective variety as a ball quotient connects with many independently interesting problems, for example the existence of special line arrangements on surfaces, the geography of Chern numbers of surfaces, and existence of extremal orbifold K\"ahler--Eienstein metrics; to only scratch the surface, see for example \cite{HirzebruchChern, HirzebruchLines, Kobayashi, DengCadorel}. Concretely, Gromov asked the following question \cite[p.\ 187]{Gromov}.

\begin{qtn}[Gromov]\label{qtn:Gromov}
Does every irreducible (smooth?) algebraic $\bbC$-variety $V$ admit a dominating regular (only rational?) map from some $\Gam \bs \bbB^n$ to some deformation $V^\prime$ of $V$?
\end{qtn}

Benson Farb asked a similar question, namely whether smooth projective surfaces are dominated by ball quotients. One perspective on the importance of Question \ref{qtn:Gromov} is that it is closely related to uniformizing varieties by orbifold metrics of constant holomorphic curvature $-1$, generalizing the classical uniformization theory for curves to higher dimensions. Indeed, essentially everything that is known about and around Question \ref{qtn:Gromov} is related to orbifold metrics modeled on the ball $\bbB^n$. Moreover, a positive answer could give insight into the classification of higher-dimensional varieties; see the discussion around Question \ref{qtn:GenType} below.

It turns out that many of the most important explicit examples of $2$-dimensional ball quotient orbifolds, constructed by Deligne and Mostow \cite{DeligneMostow}, are rational and hence birational to $\bbP^1 \times \bbP^1$. See \cite{KLW} for a precise description of the normal varieties associated with Deligne--Mostow lattices and \cite{BHH, CHL, Deraux} (and references therein) for other perspective and more on similar examples. Somewhat less well-known is Hirzebruch's construction of a compact ball quotient birational to a product of elliptic curves \cite[\S2]{HirzebruchChern}; technically, Hirzebruch considers a high degree cover and this paper extracts the simpler orbifold from his idea. This paper fully settles the rational variant of Question \ref{qtn:Gromov} for products of curves, giving the first infinite family of concrete examples for which Question \ref{qtn:Gromov} has a positive answer.

\begin{thm}\label{thm:Main}
Given $g_1, g_2 \ge 0$, there exist smooth projective curves $C_j$ of genus $g_j$ so that $C_1 \times C_2$ is birational to a compact ball quotient orbifold.
\end{thm}

It is worth noting that the birational map from the ball quotient used to prove Theorem \ref{thm:Main} to $C_1 \times C_2$ is always dominant. Smooth projective surfaces of Kodaira dimension $-\infty$ are precisely the \emph{ruled surfaces}, namely those birational to $\bbP^1 \times C$ for some smooth curve $C$ \cite[Ex.\ VII.3]{Beauville}. Therefore, Theorem \ref{thm:Main} has the following corollary.

\begin{cor}\label{cor:MainCor}
The variant of Question \ref{qtn:Gromov} where maps are birational morphisms has a positive answer for smooth projective surfaces $V$ of Kodaira dimension $\kappa(V) = -\infty$. Specifically, there is a deformation $V^\prime$ of $V$ and a ball quotient $\Gam \bs \bbB^2$ admitting a birational map $\Gam \bs \bbB^2 \dashrightarrow V^\prime$.
\end{cor}

The proof of Theorem \ref{thm:Main} consists of three primary steps. The first two are carried out in \S\ref{sec:Ex1} and \S\ref{sec:Ex2}, which construct examples for the cases $g_1 = g_2 = 2$ and $g_1 = g_2 = 1$, respectively. The ball quotient in \S \ref{sec:Ex1} is birational to $C \times C$, where $C$ is the unique smooth projective curve of genus $2$ admitting an automorphism of order $5$, first classified by Wiman \cite{Wiman}. The quotient in \S \ref{sec:Ex2} was essentially constructed by Hirzebruch, though the example in \cite[\S2]{HirzebruchChern} is technically a large degree cover of the orbifold presented in \S \ref{sec:Ex2}. The proof that each example is a ball quotient uses the orbifold Bogomolov--Miyaoka--Yau inequality, which is recalled in \S \ref{sec:Pairs}.

The third step in the proof of Theorem \ref{thm:Main} is to handle both products of higher-genus curves and ruled surfaces of genus $g \ge 2$, which are constructed by taking appropriate quotients and covers of the two primary examples. These are constructed in four stages in \S\ref{sec:Cleanup}.

Smooth projective surfaces of Kodaira dimension zero are also related to products of curves, and it was known that many can be realized as \emph{toroidal compactifications} of noncompact finite volume ball quotients. See for examples \cite{HolzapfelAbelian, KasparianKotzev, DiCerboStoverClassify, DiCerboStoverGauss}. In \S \ref{sec:Kodaira}, more new examples of ball quotient orbifolds inspired by twists on the noncompact examples prove the following result.

\begin{thm}\label{thm:Kodaira}
Let $V$ be a smooth complex projective surface of Kodaira dimension $\kappa(V) \le 0$. Then there is a compact ball quotient $\Gam \bs \bbB^n$ and a rational map $f : \Gam \bs \bbB^2 \dashrightarrow V^\prime$ to a deformation $V^\prime$ of $V$.
\end{thm}

The proofs of Theorems \ref{thm:Main} and \ref{thm:Kodaira} give even more precise results. For many cases the key step in the argument is proving that a blowup $Y$ of the appropriate $V^\prime$ is the underlying analytic space for a ball quotient orbifold, i.e., $Y$ admits a complete orbifold metric of constant biholomorphic sectional curvature $-1$. Specifically, combining this paper with previous work of Deligne--Mostow proves the existence of compact ball quotient orbifold structures on:
\begin{itemize}

\item[$\star$] del Pezzo surfaces of degree at least five

\item[$\star$] the blowup of a product $C_1 \times C_2$ of curves $C_j$ of genus $g_j \ge 2$ at a specific $3(g_1-1)(g_2-1)$ points

\item[$\star$] the blowup of an abelian surface at any finite number of points

\item[$\star$] blowups of all bielliptic surfaces except one class (the one associated with $\bbZ / 6$) at sufficiently many points

\end{itemize}
The del Pezzo cases are covered by Deligne--Mostow orbifolds \cite[Thm.\ 4.1]{KLW}. For products of higher-genus curves, see Remark \ref{rem:Higherg}. For the remaining cases, note that if $Y$ is the underlying analytic space for a ball quotient orbifold, then any \'etale covering of $Y$ is also naturally a ball quotient orbifold. The example in \S \ref{sec:Ex2} is an abelian surface blown up at one point, and the proof of Theorem \ref{thm:Kodaira} does the same for blowups of bielliptic surfaces. Appropriate \'etale covers give examples blown up at more points. The proof of Theorem \ref{thm:Kodaira} only constructs ball quotient orbifold structures on \emph{singular} K3 and Enriques surfaces, and these are the only cases in Kodaira dimension zero where this paper does not answer Question \ref{qtn:Gromov} using a dominant holomorphic map.

An example of a ball quotient orbifold structure on a smooth K3 surface was constructed in unpublished work of Naruki \cite{Naruki}; see also \cite[\S 4.4]{Hunt}. Combining Naruki's example with previously-known abelian examples, then applying the arithmeticity results of this paper described later in this introduction, proves the following result, which completely answers Farb's variant of Question \ref{qtn:Gromov} in Kodaira dimension zero.

\begin{thm}\label{thm:Farb}
Every minimal smooth projective surface of Kodaira dimension zero is smoothly dominated by a compact arithmetic ball quotient.
\end{thm}

Moreover, every simply connected $4$-manifold is dominated by the product of a torus with a higher-genus Riemann surface \cite{Branched}. The above results combine with this to give the following.

\begin{thm}\label{thm:Farb2}
Every simply connected $4$-manifold is smoothly dominated by a compact arithmetic ball quotient.
\end{thm}

As for Question \ref{qtn:Gromov} in higher Kodaira dimension, general type surfaces seem out of reach at the moment. However the Kodaira dimension one case seems quite reasonable. Surfaces of Kodaira dimension one all admit elliptic fibrations \cite[Prop.\ IX.2]{Beauville}, Theorem \ref{thm:Main} covers trivial fibrations, and Theorem \ref{thm:Kodaira} adds bielliptic surfaces. Thus the essential step in answering Question \ref{qtn:Gromov} for surfaces Kodaira dimension one is to answer the following question asked by Benson Farb.

\begin{qtn}[Farb]\label{qtn:Elliptic}
Which elliptic fibrations are dominated by a ball quotient?
\end{qtn}

For Kodaira dimension one, examples were constructed in the thesis of Livne \cite{Livne} (see also \cite[Ch.\ 16]{DMbook} and \cite{Momot}, which does not cite Livne), but little else is known beyond the trivial higher-genus fibrations constructed here. The most interesting subcase is perhaps that of Shioda modular surfaces \cite{Shioda}, which includes Livne's examples.

\medskip

Lastly, this paper studies arithmeticity of the examples used to prove the main results. Unfortunately, the examples used to prove the results in this paper do not provide new nonarithmetic lattices in $\PU(2,1)$.

\begin{thm}\label{thm:Arithmetic}
All ball quotients used to prove Theorems \ref{thm:Main}, \ref{thm:Kodaira}, and \ref{thm:Farb} are arithmetic.
\end{thm}

The example in \S \ref{sec:Ex1} is shown to be arithmetic in \S \ref{sec:Arith1}. It appears that arithmeticity of Hirzebruch's example, proved in \S \ref{sec:Arith2}, is new (Holzapfel proved arithmeticity of a different example from that paper \cite{HolzapfelPicard}). As an aside, similar techniques prove arithmeticity of a related example of Holzapfel \cite{HolzapfelAbelian} that was claimed to be arithmetic with details to appear in a later paper that does not seem to be available; see Remark \ref{rem:HolzA}. The examples from \S \ref{sec:Kodaira} are shown to be arithmetic in \S \ref{sec:LastArithmetic}. All proofs show by a geometric argument that the orbifolds are commensurable with arithmetic Deligne--Mostow orbifolds. Since all other examples in \S \ref{sec:Cleanup} are commensurable with the two primary examples, Theorem \ref{thm:Arithmetic} follows. See the end of \S \ref{sec:LastArithmetic} for the precise catalogue of results that prove Theorem \ref{thm:Arithmetic}. For a new explicit construction of a nonarithmetic orbifold (commensurable with a previously-known example), see Remarks \ref{rem:NA1} and \ref{rem:NA2}.

In fact, in light of Question \ref{qtn:Gromov}, arithmeticity can be considered a feature of this paper. A precise enumeration of smooth projective surfaces of general type seems hopeless at the current moment. However, arithmetic lattices in $\PU(2,1)$ can be enumerated and precisely classified by increasing volume, with only finitely many of a given volume bound \cite[Thm.\ A]{BorelPrasad}. Thus a positive answer to the following question could provide an intriguing conceptual step forward toward classification of surfaces of general type.

\begin{qtn}\label{qtn:GenType}
Let $V$ be a smooth projective surface of general type that is not a smooth ball quotient. Is there a deformation $V^\prime$ of $V$ and a cocompact arithmetic lattice $\Gam < \PU(2,1)$ so that $V^\prime$ is birational to the underlying analytic space for $\Gam \bs \bbB^2$?
\end{qtn}

One could also include quotients of the bidisk $\bbB^1 \times \bbB^1$ in Question \ref{qtn:GenType}, since bidisk quotients have already proven useful in constructing general type surfaces (e.g., see \cite{BCP}). Lattices acting on the bidisk are also enumerable, since they are either a product of Fuchsian groups or are arithmetic by the Margulis arithmeticity theorem \cite[p.\ 4]{MargulisBook}, so this addition would add little complication towards a goal of enumerating general type surfaces.

\subsubsection*{Acknowledgements}
I thank Benson Farb for stimulating conversations about the contents of this paper and for pointing out Naruki's K3 ball quotient \cite{Naruki} and \cite{Branched}.

\section{Ball quotient pairs}\label{sec:Pairs}

Let $\bbB^2$ be the unit ball in $\bbC^2$ with its Bergman metric and $\PU(2,1)$ its group of holomorphic isometries. A smooth projective surface $X$ is a smooth ball quotient if there is a torsion-free discrete group $\Gam < \PU(2,1)$ so that $\Gam \bs \bbB^2$ is biholomorphic to $X$. Yau famously proved that $X$ is a smooth ball quotient if and only if the canonical bundle $K_X$ is ample and the Chern numbers of $X$ satisfy $c_1^2(X) = 3 c_2(X)$, i.e., $K_X^2 = 3 e(X)$, where $e(X)$ is the topological Euler characteristic and $K_X^2$ is the self-intersection of $K_X$. Over time this result has been generalized in many ways. This paper requires a very special case of the main results in \cite{Kobayashi, KNS}.

\begin{thm}\label{thm:OrbiBMY}
Suppose $Y$ is a smooth projective surface, and $D_1, \dots, D_n$ are smooth, reduced, irreducible curves on $Y$. Given integers $r_1, \dots, r_n$, consider the $\bbQ$-divisor
\[
\calD = \sum_{j = 1}^n \left(1 - \frac{1}{r_j}\right) D_j
\]
on $Y$. Then there is a lattice $\Gam < \PU(2,1)$ so that the $\Gam \bs \bbB^2$ is an orbifold with underlying analytic space biholomorphic to $Y$ and orbifold locus $\bigcup D_j$ where $D_j$ has orbifold weight $r_j$ if and only if $K_Y + \calD$ is ample and the relative Chern numbers of the pair $(Y, \calD)$ satisfy $c_1^2(Y, \calD) = 3 c_2(Y, \calD)$.
\end{thm}

When the conditions of Theorem \ref{thm:OrbiBMY} are satisfied $(Y, \calD)$ will be called a \emph{ball quotient pair}. The relative Chern numbers in the statement of the theorem are simply
\begin{align}
c_1^2(Y, \calD) &= (K_Y + \calD)^2 \label{eq:RelativeChern1} \\
c_2(Y, \calD) &= e^{\mathrm{orb}}(Y, \calD) \label{eq:RelativeChern2}
\end{align}
where $e^{\mathrm{orb}}$ is the orbifold Euler characteristic of $Y$ with $\calD$ is considered as the orbifold locus with appropriate weights. This is subtle when curves meet, since one must know local groups at crossing points. When the curves $D_j$ are disjoint there are no crossing points with local groups to determine, hence
\begin{equation}\label{eq:Easyc2}
c_2(Y, \calD) = e(Y) - \sum_{j = 1}^n \left(1 - \frac{1}{r_j}\right) e(D_j)
\end{equation}
in that special case. The primary example in \S \ref{sec:Ex2} will have normal crossings, but is already known to be a ball quotient, and the examples in \S \ref{sec:Kodaira} will require the normal crossings variant in the exact same way: the local group at a smooth point where two curves of weights $r_j$ and $r_k$ meet transversally is simply $\bbZ / r_j \times \bbZ / r_k$.

\section{The first example}\label{sec:Ex1}

To begin, consider the smooth plane curve $C^\circ$ of genus $2$ with equation
\begin{equation}\label{eq:Curve}
y^2 = x^5 + 1
\end{equation}
in coordinates $(x,y)$ on $\bbC^2$. If $\xi$ denotes a primitive $5^{th}$ root of unity, then
\begin{equation}\label{eq:Aut}
\xi (x, y) = (\xi x, y)
\end{equation}
defines an order $5$ automorphism of $C^\circ$. The two fixed points of $\xi$ are $a_j = (0, (-1)^j)$ for $j = 1,2$. Projection $(x,y) \mapsto y$ realizes $C^\circ$ as the five-fold Galois branched cover of $\bbC$ branched over $\pm 1$ with group $\bbZ / 5 = \langle \xi \rangle$.

Let $C$ denote the unique smooth projective curve birational to $C^\circ$, which can be realized as the five-fold regular cover of $\bbP^1$ branched over $\{\pm 1, \infty\}$. The preimage of $\infty$ is a point $a_0 \in C$ so that $C \ssm \{a_0\}$ is isomorphic to $C^\circ$, and the action of $\xi$ extends to an action of $\bbZ / 5$ on $C$ with fixed points precisely at the three points $\{a_0, a_1, a_2\}$. It is originally due to Wiman in 1895 \cite{Wiman}, that $C$ is the unique smooth projective curve of genus $2$ admitting an automorphism of order $5$; see \cite{Harvey} for an alternate treatment through Fuchsian groups.

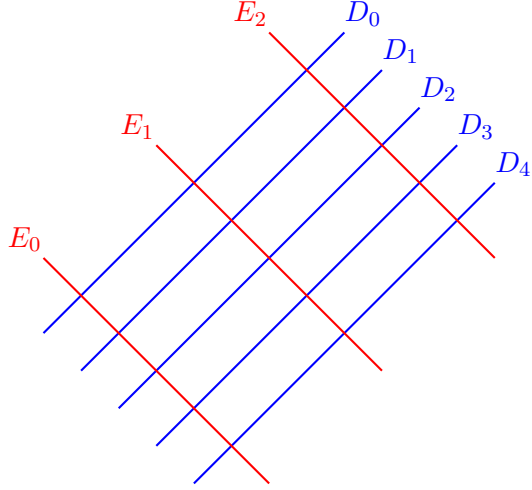
\begin{figure}[h]
\centering
\begin{tikzpicture}
\draw[thick, blue] (-2,-2) -- (2,2);
\draw[thick, blue] (-2.5,-1.5) -- (1.5,2.5);
\draw[thick, blue] (-3,-1) -- (1,3);
\draw[thick, blue] (-1.5,-2.5) -- (2.5,1.5);
\draw[thick, blue] (-1,-3) -- (3,1);
\draw[thick, red] (-1.5,1.5) -- (1.5,-1.5);
\draw[thick, red] (-3,0) -- (0,-3);
\draw[thick, red] (3,0) -- (0,3);
\node[thick, red] at (-1.75, 1.75) {$E_1$};
\node[thick, red] at (-3.25, 0.25) {$E_0$};
\node[thick, red] at (-0.25, 3.25) {$E_2$};
\node[thick, blue] at (2.25, 2.25) {$D_2$};
\node[thick, blue] at (1.75, 2.75) {$D_1$};
\node[thick, blue] at (1.25, 3.25) {$D_0$};
\node[thick, blue] at (2.75, 1.75) {$D_3$};
\node[thick, blue] at (3.25, 1.25) {$D_4$};
\end{tikzpicture}
\caption{The arrangement of curves on $Y$}\label{fig:Arrangement}
\end{figure}

Set $X = C \times C$, and let $C_j$ be the graph of $\xi^j$ for $0 \le j \le 4$. Since $\xi^{k-j}$ has the same fixed points as $\xi$ for any $j \neq k$,
\[
C_j \cap C_k = \{x_0, x_1, x_2\}
\]
where $x_j = (a_j, a_j)$. Let $\pi : Y \to X$ be the blowup of $X$ at $\{x_0, x_1, x_2\}$, $E_j$ denote the exceptional divisor associated with $x_j$, and $D_j$ be the proper transform of $C_j$ to $Y$. See Figure \ref{fig:Arrangement} for a depiction of this simple arrangement of curves on $Y$ and \cite[\S II.1]{Beauville} for basic properties of blowups that will be used in what follows.

Let $C_h$ and $C_v$ denote a general horizontal and vertical divisor of $X$, respectively. The canonical divisor $K_X$ of $X$ is then $2(C_h + C_v)$ by \cite[II Ex.\ 8.3]{Hartshorne} and the fact that the canonical divisor on a smooth projective curve of genus $g$ has degree $2 g - 2$. If $D_h$ and $D_v$ are the proper transforms of $C_h$ and $C_v$ to $Y$, then
\begin{equation}\label{eq:Ycanonical}
K_Y = 2(D_h + D_v) + \sum E_j.
\end{equation}
Note that $K_Y^2 = 5$. Since $C_j \cdot C_h = C_j \cdot C_v = 1$ for all $j$ and $D_j \cdot E_k = 1$ for all $j,k$, it follows that $K_Y \cdot D_j = 7$, and thus $D_j^2 = -5$ by adjunction \cite[I.15]{Beauville}. These calculations are now used to prove the following result, which provides the first of the two examples bootstrapped to prove Theorem \ref{thm:Main}.

\begin{thm}\label{thm:Ex}
Let $Y$ be the surface described previously in this section, $\{D_j\}$ be the curves of genus $2$ defined on $Y$, and consider the $\bbQ$-divisor
\begin{equation}\label{eq:Divisor}
\calD = \sum_{j = 0}^4 \left(1 - \frac{1}{5}\right) D_j
\end{equation}
on $Y$. Then $(Y, \calD)$ is a ball quotient pair.
\end{thm}

\begin{pf}
Define $\calL = K_Y + \calD$. To apply Theorem \ref{thm:OrbiBMY}, the first task is to compute the Chern numbers $c_1^2(Y, \calD)$ and $c_2(Y, \calD)$. From earlier in this section, $K_Y^2 = 5$, $K_Y \cdot D_j = 7$, and $D_j^2 = -5$ for each $0 \le j \le 4$, hence it follows that
\begin{equation}\label{eq:c1}
c_1^2(Y, \calD) = 45.
\end{equation}
Since the curves $D_j$ are disjoint of genus $2$ and $c_2(Y) = 7$, Equation \eqref{eq:Easyc2} gives
\begin{equation}\label{eq:c2}
c_2(Y) - \sum_{j = 0}^4 \left(1 - \frac{1}{5} \right) (-2) = 15.
\end{equation}
Thus $c_1^2(Y, \calD) = 3 c_2(Y, \calD)$ as required.

It remains to show that $\calL$ is ample. By the Nakai--Moishezon criterion \cite[Thm.\ V.1.10]{Hartshorne}, since $\calL^2 > 0$ it suffices to show that $\calL \cdot A > 0$ for every irreducible curve $A$ on $Y$. As $Y$ is the blowup of the minimal surface $X$ of general type, $K_X$ is ample and so $K_Y \cdot A > 0$ for any irreducible curve $A$ on $Y$ that is not one of the three exceptional curves. Since distinct curves have nonnegative intersection number, $\calL \cdot A > 0$ when $A$ is not one of the three exceptional curves or one of the curves $D_j$. It therefore suffices to show that $\calL \cdot A > 0$ when $A$ is either an exceptional curve or one of the curves $D_j$. Direct computation gives
\[
\calL \cdot E_j = \calL \cdot D_j = 3
\]
for all relevant $j$, and this completes the proof.
\end{pf}

\section{Arithmeticity of the example in \S \ref{sec:Ex1}}\label{sec:Arith1}

Let $(Y, \calD)$ be the ball quotient pair in Theorem \ref{thm:Ex} and $\Gam \bs \bbB^2$ be the associated ball quotient orbifold. Continuing with the notation of \S \ref{sec:Ex1}, consider the action of $G = (\bbZ / 5)^2$ on $X = C \times C$ given by $(z, w) \mapsto (\xi^r z, \xi^s w)$. Since the points $x_j = (a_j, a_j) \in X$ are global fixed points of the action, where $a_j \in \{\pm 1, \infty\}$ are the fixed points of $\xi$, the action of $G$ lifts to an action on the blowup $Y$ stabilizing the three exceptional curves $E_j$. Let $Z$ be the quotient $Y / G$.

\begin{lem}\label{lem:ZisP2}
The surface $Z$ is isomorphic to $\bbP^2$ blown up in the vertices of the complete quadrangle.
\end{lem}

\begin{pf}
Recall that the complete quadrangle in $\bbP^2$ consists of the six lines where either two coordinates are equal or one is equal to zero; see \cite[Fig.\ I(vi)]{KLW} or Figure \ref{fig:P123} below for the standard depiction of the quadrangle. Since $C / \langle \xi \rangle$ is isomorphic to $\bbP^1$, $X / G$ is isomorphic to $\bbP^1 \times \bbP^1$. For each point $y$ on an exceptional curve $E_j$, one checks directly in coordinates that its stabilizer $\Stab_G(y)$ is generated by complex reflections, hence $Y / G$ is a smooth surface by the Chevalley--Shephard--Todd theorem \cite[\S V.5.5, Th.\ 4]{BourbakiV} and $E_j$ maps to a smooth curve $\wh{E}_j$ on $Z$ isomorphic to $\bbP^1$. Moreover, contraction of $\{\wh{E}_j\}$ is a map $Z \to \bbP^1 \times \bbP^1$, and therefore $\wh{E}_j$ is a $(-1)$-curve; see \cite[\S III.4]{BPV}. This proves that $Z$ is isomorphic to $\bbP^1 \times \bbP^1$ blown up at three points on the diagonal, and it is a classical exercise that this surface is isomorphic to $\bbP^2$ blown up at the four vertices of the complete quadrangle.
\end{pf}

Moreover, there is a $\bbQ$-divisor $\calE$ on $Z$ so that $(Z, \calE)$ is a ball quotient pair and the associated ball quotient orbifold $\Lam \bs \bbB^2$ is isomorphic to $(\Gam \bs \bbB^2) / G$. To prove this, let $\Aut_\calD(Y)$ denote the automorphisms of $Y$ preserving $\calD$, in that sense that it not only preserves the support of $\calD$ but also preserves the weights of divisors. By the definition of an orbifold cover, the quotient $Y / \Aut_\calD(Y)$ is again a ball quotient and the induced map is an orbifold covering.

\begin{lem}\label{lem:Gquotient}
The action of $G$ on $Y$ induces an action on $\Gam \bs \bbB^2$, hence the quotient $(\Gam \bs \bbB^2) / G$ is a ball quotient $\Lam \bs \bbB^2$ for a lattice $\Lam < \PU(2,1)$ containing $\Gam$ as a normal subgroup with quotient $G$. The orbifold $\Lam \bs \bbB^2$ arises from a ball quotient pair $(Z, \calE)$ for some $\bbQ$-divisor $\calE$ on $Z$.
\end{lem}

\begin{pf}
Recall that the curves $D_j$ on $Y$ are the proper transforms of the graphs of $\xi^j$, hence are isomorphic to $C$. The action of $G$ cyclically permutes these curves, and the stabilizer of each $D_j$ is isomorphic to $\bbZ / 5$ acting on $C$ by the usual action. This alone suffices to imply that $G$ acts on $\Gam \bs \bbB^2$, since the $G$-action on $Y$ preserves the orbifold locus. All statements in the lemma are restatements of this fact and hence follow directly.
\end{pf}

The next goal is to explicitly determine the $\bbQ$-divisor $\calE$ on $Z$ and hence (equivalently) describe the orbifold structure on $Z$ as $\Lam \bs \bbB^2$. The diagonal subgroup of $G$ acts on an exceptional curve $E_j$ as a complex reflection through $E_j$ and the quotient of $G$ by the diagonal subgroup acts on $E_j$ as $y \mapsto y^5$. The remaining curves stabilized by $G$ are the proper transforms $H_j, V_j$ ($j \in \{0,1,2\}$) of the horizontal and vertical curves on $X$ with fixed coordinate $a_j$. The action of $G$ on any of these curves $H_j, V_j$ is generated by a complex reflection of order $5$ and the usual action of $\bbZ / 5$ on $C$.

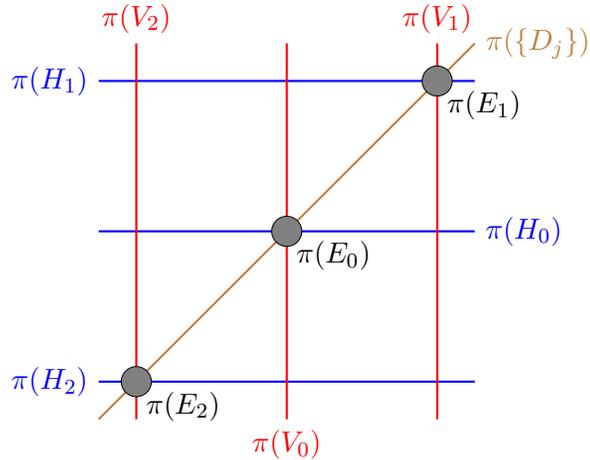
\begin{figure}[h]
\centering
\begin{tikzpicture}
\draw[thick, blue] (-2.5,2) node[left] {$\pi(H_1)$} -- (2.5,2);
\draw[thick, blue] (-2.5,0) -- (2.5,0) node[right] {$\pi(H_0)$};
\draw[thick, blue] (-2.5,-2) node[left] {$\pi(H_2)$} -- (2.5,-2);
\draw[thick, red] (-2,-2.5) -- (-2,2.5) node[above] {$\pi(V_2)$};
\draw[thick, red] (0,-2.5) node[below] {$\pi(V_0)$} -- (0,2.5) ;
\draw[thick, red] (2,-2.5) -- (2,2.5) node[above] {$\pi(V_1)$};
\draw[thick, brown] (-2.5, -2.5) -- (2.5, 2.5) node [right] {$\pi(\{D_j\})$};
\draw[very thin, fill = gray] (-2,-2) node[below right] {$\pi(E_2)$} circle (0.2cm);
\draw[very thin, fill = gray] (0,0) node[below right] {$\pi(E_0)$} circle (0.2cm);
\draw[very thin, fill = gray] (2,2) node[below right] {$\pi(E_1)$} circle (0.2cm);
\end{tikzpicture}
\caption{The arrangement of curves on $Z$}\label{fig:QuoPic}
\end{figure}

See Figure \ref{fig:QuoPic} for a depiction of the arrangement of curves on $Z$ defined by the images of the curves $D_j$, $E_j$, $H_j$, and $V_j$. The branch divisor for the quotient map $\pi$ is the union of the $E_j$, $H_j$, and $V_j$, and the image of the curves $D_j$ is the proper transform to $Z$ of the diagonal of $\bbP^1 \times \bbP^1$ with respect to the natural blowup map $Z \to \bbP^1 \times \bbP^1$ commuting with the action of $G$.

To understand the ball quotient orbifold structure on $Z$, first recall that the orbifold structure on $Y$ consists of the five curves $D_j$, each with orbifold weight $5$. The curves $D_j$ are not in the branch locus of $\pi$, so the image on $Z$ has orbifold weight $5$ on $\Lam \bs \bbB^2$. The remaining curves in the arrangement are not in the orbifold locus on $Y$, but $\pi$ branches to order $5$ along each, so it becomes a curve of orbifold weight $5$ in the ball quotient structure on $Z$. In other words, the orbifold locus on $Z$ consists of precisely the curves in the arrangement in Figure \ref{fig:QuoPic}, each with weight $5$. Said yet one more way, if
\[
\calE = \left(1 - \frac{1}{5}\right) \left( \pi(\{D_j\}) + \sum_{j = 0}^2 \left(\pi(E_j) + \pi(H_j) + \pi(V_j)\right)\right)
\]
then $(Z, \calE)$ is the ball quotient pair associated with $\Lam \bs \bbB^2$ and the action of $G$ on $Y$ induces an orbifold covering $\Gam \bs \bbB^2 \to \Lam \bs \bbB^2$.

Using the standard identification of the arrangement in Figure \ref{fig:QuoPic} with the blowup of the complete quadrangle in $\bbP^2$ at its vertices, this implies that $(Z, \calE)$ can be realized as the arrangement shown in \cite[Fig.\ I(i)]{KLW}. Factoring in the orbifold weight of each curve, this is the Deligne--Mostow quotient with weights $(2,2,2,2,2)/5$ by \cite[Thm.\ 4.1]{KLW}, which is arithmetic \cite[p.\ 86]{DeligneMostow}. Note that $(Z, \calE)$ is associated with the quotient of $\bbB^2$ by the `pure' Deligne--Mostow lattice, not its extension by the action of the symmetric group $S_6$ on the weights. This completes the proof of the following result.

\begin{thm}\label{thm:Ex1A}
The ball quotient orbifold constructed in \S \ref{sec:Ex1} is arithmetic.
\end{thm}

\section{The second example}\label{sec:Ex2}

Let $T$ be the elliptic curve $\bbZ[\zeta] \bs \bbC$, where $\zeta = e^{2 \pi i / 6}$, and consider the abelian surface $T \times T$. In coordinates $(z, w)$ on $T \times T$, let $T_j$ be the graph of multiplication by $j \in \bbZ[\zeta]$. As noted in \cite[\S 1]{HirzebruchChern}, the curves $T_0$, $T_\infty$ (i.e., $z = 0$), $T_1$, and $T_\zeta$ are embedded elliptic curves that mutually intersect at the origin $(0,0)$ and otherwise in no other point.

Let $Y$ be the blowup of $T \times T$ at $(0,0)$, $L$ be the exceptional curve, $D_j$ denote the proper transform of $T_j$ to $Y$, and define
\[
\calD =\! \left(1 - \frac{1}{3}\right)\!\left(L + \sum D_j\right)\!.
\]
The following proposition is essentially equivalent to the result proved in \cite[\S 2]{HirzebruchChern}, and hence covers the case $g_1 = g_2 = 1$ in Theorem \ref{thm:Main}. A proof based on \S \ref{sec:Pairs} is provided for completeness.

\begin{prop}\label{prop:HirzebruchPair}
The pair $(Y, \calD)$ defined in this section is a ball quotient pair.
\end{prop}

\begin{pf}
Note that $K_Y$ is numerically equivalent to the exceptional curve $L$ (e.g., see \cite[\S II.1]{Beauville}). Since the $D_j$ are disjoint, each $D_j$ has self-intersection $-1$, and $L \cdot D_j = 1$ for all $j$ it follows that:
\begin{align*}
c_1^2(Y, \calD) &= (L + \calD)^2 \\
&=\!\left(\frac{5}{3} L + \frac{2}{3} \sum D_j\right)^2 \\
&= \frac{13}{3} \\
c_2(Y, \calD) &= e(Y) - \left(1 - \frac{1}{3}\right) e\!\left(L \ssm \{4 \textrm{pts}\}\right) \\
& \quad - \sum \left(1 - \frac{1}{3}\right) e\!\left(D_j \ssm \{1 \textrm{pt}\}\right) - \left(1 - \frac{1}{9}\right) e\!\left(\{4 \textrm{pts}\}\right) \\
&= 1 - \frac{2}{3}(-2) - \frac{8}{3}(-1) - \frac{32}{9} \\
&= \frac{13}{9}
\end{align*}
Thus $c_1^2(Y, \calD) = 3 c_2(Y, \calD)$.

It remains to prove that $\calL = K_Y + \calD$ is ample, which follows the same script as the proof of Theorem \ref{thm:Ex}. Notice that $\calL$ is an effective linear combination of $L$ and the $D_j$. Any curve on $Y$ not $L$ or some $D_j$ meets some curve in the support of $\calL$. Indeed, suppose that $C$ does not meet $D_0$ or $D_\infty$. Then the image of $C$ in $T \times T$ must meet either $C_0$ or $C_\infty$, which means that $C$ meets $L$. It follows that any such $C$ meets some curve in the support of $\calD$, and thus $\calL \cdot C > 0$. Finally, the computations
\begin{align*}
L\cdot \calL &= \frac{5}{3}(-1) + \frac{8}{3} = 1 \\
D_j \cdot \calL &= \frac{5}{3} + \frac{2}{3}(-1) = 1
\end{align*}
then suffice by the Nakai--Moishezon criterion \cite[Thm.\ V.1.10]{Hartshorne} to imply that $\calL$ is ample. Thus $(Y, \calD)$ is a ball quotient pair.
\end{pf}

\begin{rem}
The manifold $W_3$ in \cite[\S 2]{HirzebruchChern} is a degree $3^8$ cover of the orbifold given by Proposition \ref{prop:HirzebruchPair}.
\end{rem}

\begin{rem}
The examples provided by Theorem \ref{thm:Ex} and Proposition \ref{prop:HirzebruchPair} are closely related. Indeed, both are constructed using the unique smooth projective curve of genus $g \in \{2, 1\}$ admitting an automorphism of order $2 g + 1$; in fact, both admit one of order $4 g + 2$, and they are the unique curves of their genus achieving this maximum. A systematic computer search for configurations to which Theorem \ref{thm:OrbiBMY} applies found no obvious analogous higher genus ball quotient pair.
\end{rem}

\section{Arithmeticity of the example in \S \ref{sec:Ex2}}\label{sec:Arith2}

Let $(Y, \calD)$ be as in \S \ref{sec:Ex2}, $\Gam \bs \bbB^2$ be the associated ball quotient, and retain all other notation from that section. As in \S \ref{sec:Arith1}, $\Gam$ is shown to be arithmetic is by proving that a certain quotient of $\Gam \bs \bbB^2$ is an arithmetic Deligne--Mostow orbifold. The argument this time is somewhat more complicated, since the quotient is a singular rational surface that is not as simple to identify. While there may be an identification of the quotient relying on more classical techniques, this section uses an auxiliary ball quotient, the noncompact manifold constructed by Hirzebruch in the same paper \cite[p.\ 356]{HirzebruchChern}, which admits an action by the same finite group. This argument has the added bonus of giving a natural explanation for a certain group of order $72$, previously found only computationally, relating Hirzebruch's noncompact example to a Picard modular surface that is also a Deligne--Mostow orbifold \cite[Appendix]{StoverUrzua}.

Let $\mu_6 = \langle \zeta \rangle$ be the group of $6^{th}$ roots of unity, and $F$ be the subgroup of $\GL_2(\mu_6)$ generated by:
\begin{align*}
\al &= \begin{pmatrix} 1 & \zeta^2 \\ 1 & -1 \end{pmatrix} &
\beta &= \begin{pmatrix} 1 & 0 \\ 1 & \zeta^2 \end{pmatrix}
\end{align*}
There is a natural action of $F$ on $T \times T$ fixing the origin, and hence an action on the blowup $Y$. The action of $\al$ on the curves $T_j$ is:
\begin{align*}
\al(z, 0) &= (z, z) & \al(z, z) &= (-\zeta z, 0) \\
\al(0, z) &= (\zeta^2 z, -z) & \al(z, \zeta z) &= (0, -\zeta^2 z)
\end{align*}
Noting that $(\zeta^2 z, -z) = (w, \zeta w)$ for $w = \zeta^2 z$, $\al$ permutes the $T_j$ by the permutation $(1\, 3)(2\, 4)$ for the ordering $\{T_0, T_\infty, T_1, T_\zeta\}$. Moreover, $\al$ has order $12$ with
\[
\al^2 = \begin{pmatrix} \zeta & 0 \\ 0 & \zeta \end{pmatrix} = \zeta \Id
\]
central in $F$. Further calculations show that $\beta$ acts on the $T_j$ as $(1\, 3\, 4)$, $\beta^3$ is the identity on $T \times T$, $(\al \beta)$ acts as $(2\, 4\, 3)$ on the $\{T_j\}$, and $(\al \beta)^3 = -\Id$ on $T \times T$. Therefore the restriction to permutations of the $\{T_j\}$ defines a homomorphism of $F$ onto the alternating group $A_4$ with kernel $\zeta \Id$.

\begin{rem}
This will not be used in what follows, but one can deduce an abstract presentation
\[
F =\! \left\langle \al, \beta\ |\ \al^{12}, \beta^3, [\al^2, \beta], (\al \beta)^3 = \al^6 \right\rangle
\]
with $\langle \al^2 \rangle \cong \bbZ / 6$ central and $F / \langle \al^2 \rangle \cong A_4$. In other words, $F$ is a central extension of $A_4$ by $\bbZ / 6$.
\end{rem}

The induced action of $F$ on $Y$ also permutes the proper transforms $D_j$ of the curves $T_j$ and the action on the exceptional curve $L$ is by the linear action of the given matrix on lines through the origin of $T \times T$. Exactly as in Lemma \ref{lem:Gquotient}, the action of $F$ on $Y$ preserves the orbifold structure for the pair $(Y, \calD)$, hence $F \le \Aut_\calD(Y)$ and the following lemma holds.

\begin{lem}\label{lem:Fquotient}
The action of $F$ on $Y$ induces an action on $\Gam \bs \bbB^2$. Therefore there is a lattice $\Lam < \PU(2,1)$ containing $\Gam$ as a normal subgroup of index $|F| = 72$ so that $\Lam \bs \bbB^2 = (\Gam \bs \bbB^2) / F$.
\end{lem}

The goal is now to determine the space $Y / F$ and prove that $\Lam$ is an arithmetic Deligne--Mostow lattice. The first step is to show that the quotient is isomorphic to a scheme-theoretic blowup of the weighted projective plane $\bbP(1,2,3)$ by first proving that $(Y \ssm \bigcup D_j) / F$ is isomorphic to the complement of a certain point in $\bbP(1,2,3)$. This open variety is the underlying space for a \emph{different} Deligne--Mostow orbifold.

Consider the divisor $\calD_\infty = \sum D_j$ on $Y$, which fits into \S \ref{sec:Pairs} in the more general setting of weight $r_j = \infty$ that allows for cusps. It is now well-established in the literature that $(Y, \calD_\infty)$ is a ball quotient pair in this sense, i.e., that $Y^\circ = Y \ssm \bigcup D_j$ is the quotient $\Del \bs \bbB^2$ of $\bbB^2$ by an arithmetic lattice $\Del < \PU(2,1)$. For example, see Example 1 in \S 6 and \S 7 of \cite{DiCerboStoverClassify}. The exact same logic that was applied to $Y$ and $\Gam \bs \bbB^2$ gives the following lemma.

\begin{lem}\label{lem:Fquotient2}
The action of $F$ on $Y^\circ$ induces an action on $\Del \bs \bbB^2$. Therefore there is a lattice $\wh{\Del} < \PU(2,1)$ containing $\Del$ as a normal subgroup of index $|F| = 72$ so that $\wh{\Del} \bs \bbB^2 = (\Del \bs \bbB^2) / F$.
\end{lem}

The identification of $Y^\circ / F$ with $\bbP(1,2,3)$ minus a point will come from the following proposition, which identifies the group $\wh{\Del}$.

\begin{prop}\label{prop:Delhat}
The group $\wh{\Del}$ in Lemma \ref{lem:Fquotient2} is the Deligne--Mostow lattice with weights $(5,4,1,1,1)/6$.
\end{prop}

\begin{pf}
It is known that $\Del$ is commensurable with the Picard modular group $\PU(2,1; \bbZ[\zeta])$ \cite[\S 7]{DiCerboStoverClassify}. A direct calculation shows that $\Del \bs \bbB^2$ has Euler characteristic $1$, hence $\wh{\Del} \bs \bbB^2$ has orbifold Euler characteristic $1/72$, which is the smallest possible orbifold Euler characteristic of a noncompact arithmetic quotient of $\bbB^2$ by \cite[Thm.\ 1.1]{StoverVols} and Chern--Gauss--Bonnet. It is already known that $\Gam$ is a normal subgroup with index $72$ in the Deligne--Mostow lattice with weights $(5,4,1,1,1)/6$ \cite[Appendix]{StoverUrzua}, and this must be $\wh{\Del}$ by Mostow rigidity.
\end{pf}

\begin{rem}
The weights $(5,4,1,1,1)/6$ only satisfy condition $\Sig\mathrm{INT}$ of Deligne and Mostow, hence there is no ambiguity of `pure' lattice versus the extension by the $S_3$ action on weights.
\end{rem}

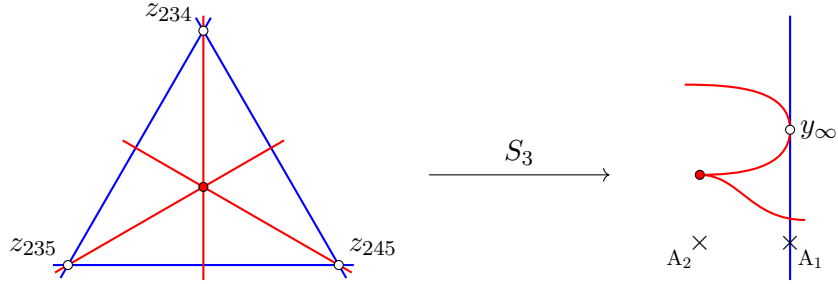
\begin{figure}[h]
\centering
\begin{tikzpicture}[scale=1.2]
\draw[blue, thick, shorten >= -0.2cm, shorten <= -0.2cm] (-5, 0) -- (-2, 0);
\draw[blue, thick, shorten >= -0.2cm, shorten <= -0.2cm] (-5, 0) -- (-3.5, 2.598);
\draw[blue, thick, shorten >= -0.2cm, shorten <= -0.2cm] (-2, 0) -- (-3.5, 2.598);
\draw[red, thick, shorten >= -0.2cm, shorten <= -0.2cm] (-3.5, 2.598) -- (-3.5, 0);
\draw[red, thick, shorten >= -0.2cm, shorten <= -0.2cm] (-4.25, 1.299) -- (-2, 0);
\draw[red, thick, shorten >= -0.2cm, shorten <= -0.2cm] (-2.75, 1.299) -- (-5, 0);
\draw[fill=white] (-5,0) circle (0.05cm) node [above left] {$z_{235}$};
\draw[fill=white] (-2,0) circle (0.05cm) node [above right] {$z_{245}$};
\draw[fill=white] (-3.5, 2.598) circle (0.05cm) node [above left] {$z_{234}$};
\draw[fill=red] (-3.5, 0.866) circle (0.05cm);
\draw[->] (-1,1) -- node [above] {$S_3$} (1,1);
\draw[blue, thick, shorten >= -0.2cm, shorten <= -0.2cm] (3, 0) -- (3, 2.598);
\draw [red, thick, smooth, shorten >= -0.2cm, shorten <= -0.2cm] (2,2) to[out=0, in=90] (3,1.5) to[out=270, in=0] (2,1) to[out=0, in=180] (3,0.5);
\draw[fill=white] (3,1.5) circle (0.05cm) node [right] {$y_\infty$};
\draw[fill=red] (2,1) circle (0.05cm);
\draw (2,0.25) node {$\times$} node [below left, scale=0.75] {$\mathrm{A}_2$};
\draw (3,0.25) node {$\times$} node [below right, scale=0.75] {$\mathrm{A}_1$};
\end{tikzpicture}
\caption{The map $(\bbP^2 \ssm \{z_{234}, z_{235}, z_{245}\}) \to (\bbP(1,2,3) \ssm \{y_\infty\})$}\label{fig:P123}
\end{figure}

\begin{cor}\label{cor:P123}
The quotient $Y^\circ / F$ is the complement of the point $y_\infty$ in $\bbP(1,2,3)$ indicated in Figure \ref{fig:P123}.
\end{cor}

\begin{pf}
Following Theorem 4.1(vi), (5.1)(vi), and Figure I(vi) in \cite{KLW}, $\wh{\Del} \bs \bbB^2$ has underlying space the quotient $(\bbP^2 \ssm \{z_{234}, z_{235}, z_{245}\}) / S_3$, where $z_{ijk}$ is the vertex of the complete quadrangle in $\bbP^2$ indicated in Figure \ref{fig:P123} and $S_3$ acts by coordinate permutations. The quotient of $\bbP^2$ by $S_3$ is classically known to be $\bbP(1,2,3)$ (e.g., see \cite[Ch.\ 11]{DMbook}), where the image of the complete quadrangle and the singularities of type $\mathrm{A}_1$ and $\mathrm{A}_2$ are depicted in Figure \ref{fig:P123}. This proves the lemma.
\end{pf}

The above information is now used to prove the following result.

\begin{prop}\label{prop:FoundDM}
The space $\Lam \bs \bbB^2$ in Lemma \ref{lem:Fquotient} is the Deligne--Mostow quotient with weights $(14,13,3,3,3)/18$, which has underlying analytic space the scheme theoretic blowup of $\bbP(1,2,3)$ at $y_\infty$.
\end{prop}

\begin{pf}
The complement in $\Gam \bs \bbB^2$ of the four curves $D_j$ in the orbifold locus is isomorphic to $Y^\circ$ as a complex manifold. The action of $F$ on the $D_j$ is transitive, so they map to a single curve $\wh{D}$ in the orbifold $\Lam \bs \bbB^2$ in Lemma \ref{lem:Fquotient}. By Corollary \ref{cor:P123}, the quotient
\[
\left(\Gam \bs \bbB^2 \ssm \bigcup D_j\right)\!/F
\]
is isomorphic to $\bbP(1,2,3) \ssm \{y_\infty\}$, and hence $\Lam \bs \bbB^2$ has underlying space obtained by replacing $y_\infty$ with the curve $\wh{D}$. Thus must be made precise.

The stabilizer $F_j$ in $F$ of a curve $D_j$ is isomorphic to $\bbZ / 6 \times \bbZ / 3$ generated by the center $\langle \zeta \Id \rangle$ along with a lift of the appropriate $3$-cycle in $\mathrm{A}_4$. That $3$-cycle acts on $D_j$ by $\zeta^2$, so there is a $\bbZ / 3$ subgroup $R_j < F_j$ that acts trivially on $D_j$, and $F_j / R_j$ is the usual action of $\zeta$ on $D_j$. In other words, the image of $D_j$ on $\Lam \bs \bbB^2$ is isomorphic to $T / \langle \zeta \rangle \cong \bbP^1$. The orbifold weight of $D_j$ on $\Gam \bs \bbB^2$ is $3$, and $F_j$ has subgroup of order $3$ acting as a complex reflection through $D_j$, so $\wh{D}$ has orbifold weight $9$ on $\Lam \bs \bbB^2$. 

The only fixed point of $\zeta$ on $T$ is the origin, which on $D_j$ is its point of intersection with the exceptional curve $L$. The action of $\zeta$ on $L$ is trivial, hence $\zeta$ locally acts as a complex reflection at $D_j \cap L$, and thus this fixed point projects to a smooth point of $\Lam \bs \bbB^2$ by the Chevalley--Shephard--Todd theorem \cite[\S V.5.5, Th.\ 4]{BourbakiV}. The action of $-1 = \zeta^3$ on $D_j$ has three additional fixed points at nontrivial $2$-torsion elements of $T$, where the action on $Y$ is locally by $-\Id$. These three points are permuted by the action of $\zeta^2$, hence they determine a single $\mathrm{A}_1$ singularity on $Y/F$ that lies on $\wh{D}$. The two additional fixed points of $\zeta^2$, which are nontrivial $3$-torsion points identified by the action of $\zeta^3$, have stabilizer $(\bbZ / 3)^2 < F_j$ generated by complex reflections, hence they descend to a smooth point of $Y / F$.

In conclusion, $\wh{D}$ is isomorphic to $\bbP^1$, passes through an $\mathrm{A}_1$ singularity of $Y / F$, and has orbifold weight $9$ for the ball quotient structure on $\Lam \bs \bbB^2$. Contraction of $\wh{D}$ is a map to $\bbP(1,2,3)$. The universal property of blowing up \cite[Prop.\ II.7.14]{Hartshorne} implies that this map factors through the scheme theoretic blowup of $\bbP(1,2,3)$ at $y_\infty$, but $\wh{D}$ is reduced and irreducible, hence $Y / F$ must be isomorphic to the blowup. See Figure \ref{fig:Y/F} and compare with ${}_2 R$ in the figure on page 113 of \cite{DMbook}. It remains to determine the rest of the orbifold locus of $\Lam \bs \bbB^2$, i.e., the image on $Y / F$ of the orbifold locus of $\Gam \bs \bbB^2$ and any additional ramification of the projection.
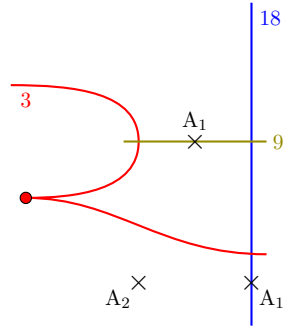
\begin{figure}[h]
\centering
\begin{tikzpicture}[scale=1.5]
\draw[blue, thick, shorten >= -0.2cm, shorten <= -0.2cm] (3, 0) -- (3, 2.598) node [right, scale=0.75] {$18$};
\draw[red, thick, smooth, shorten >= -0.2cm, shorten <= -0.2cm] (1,2) node [below, scale=0.75] {$3$} to[out=0, in=90] (2,1.5) to[out=270, in=0] (1,1) to[out=0, in=180] (3,0.5);
\draw[olive, thick, shorten >= -0.2cm, shorten <= -0.2cm] (2,1.5) -- (3,1.5) node [right, xshift=5, scale=0.75] {$9$};
\draw[fill=red] (1,1) circle (0.05cm);
\draw (2,0.25) node {$\times$} node [below left, scale=0.75] {$\mathrm{A}_2$};
\draw (3,0.25) node {$\times$} node [below right, scale=0.75] {$\mathrm{A}_1$};
\draw (2.5,1.5) node {$\times$} node [above, yshift=1.5, scale=0.75] {$\mathrm{A}_1$};
\end{tikzpicture}
\caption{The quotient $Y / F$ with orbifold weights}\label{fig:Y/F}
\end{figure}

The only other curve in the orbifold locus of $\Gam \bs \bbB^2$ is the exceptional curve $L$, which has orbifold weight $3$ and is stabilized by $F$. The central $\langle \zeta \Id \rangle < F$ acts by a complex reflection through $L$, since scalar matrices act trivial on lines through the origin. The action of $F / \langle \zeta \Id \rangle \cong A_4$ is by the standard action of $A_4$ on $\bbP^1$. Thus $L$ projects to a smooth rational curve on $L / F$ with orbifold weight $18$. One can see either through analysis of the action of $F$ or by noticing that $L$ maps to the curve of orbifold weight $6$ for the interpretation of $Y^\circ / F$ as the $(5,4,1,1,1)/6$ Deligne--Mostow orbifold that the image of $L$ is the curve in Figure \ref{fig:Y/F} passing through the second $\mathrm{A}_1$ singularity of $Y / F$.

Finally, any additional ramification of $Y \to Y / F$ arises from complex reflections in $F$. Every element of $F$ is conjugate to $\zeta^r \al$ or $\zeta^r \beta$ for some $r$, thus the only complex reflections in $F$ are the conjugates of $\beta$. Computing the eigenvectors of $\beta$, the curve $T_\beta$ on $T \times T$ fixed by $\beta$ is the embedded $T$ with parametrization
\[
z \longmapsto\!\left((1-\zeta^2) z, z\right)\!.
\]
The full stabilizer of $T_\beta$ in $F$ also contains the center, which again acts on $T_\beta$ as multiplication by $\zeta$. Either similar analysis as with $\wh{D}$ and $L$ or analogy with the $(5,4,1,1,1)/6$ Deligne--Mostow orbifold implies that the image of $T_\beta$ on $Y / F$ is the singular curve in Figure \ref{fig:Y/F} with orbifold weight $3$ on $\Lam \bs \bbB^2$.

Therefore Figure \ref{fig:Y/F} depicts all the orbifold locus on $\Lam \bs \bbB^2$. Following \cite[Ch.\ 11]{DMbook} and \cite[Thm.\ 4.1(ii)]{KLW}, this is precisely the Deligne--Mostow orbifold with weights $(14,13,3,3,3)/18$. This completes the proof of the proposition.
\end{pf}

Since the Deligne--Mostow orbifold with weights $(14,13,3,3,3)/18$ is arithmetic \cite[p.\ 102]{MostowINT}, the main result of this section now follows immediately from Proposition \ref{prop:FoundDM}.

\begin{thm}\label{thm:Ex2A}
The ball quotient orbifold constructed in \S \ref{sec:Ex2} is arithmetic.
\end{thm}

\section{Remaining cases for Theorem \ref{thm:Main}}\label{sec:Cleanup}

The results in \S \ref{sec:Ex1} and \S \ref{sec:Ex2} provide the cases $g_1 = g_2 = 2$ and $g_1 = g_2 = 1$ in the proof of Theorem \ref{thm:Main}, respectively. The case $g_1 = g_2 = 0$ is handled by Deligne--Mostow orbifolds (e.g., see the proof of Theorem \ref{thm:Ex1A}). This section bootstraps from those examples to complete the proof of Theorem \ref{thm:Main}.

\begin{prop}\label{prop:MainHyp}
Given $g_1, g_2 \ge 2$, there exist smooth projective curves $C_j$ of genus $g_j$ so that $C_1 \times C_2$ is birational to a ball quotient.
\end{prop}

\begin{pf}
Let $\Gam < \PU(2,1)$ be the lattice determined by the example in Theorem \ref{thm:Ex} and retain all notation from \S\ref{sec:Ex1}. Forgetting the orbifold structure on $\Gam \bs \bbB^2$ induces an exact sequence
\[
1 \lra R \lra \Gam \lra \pi_1(Y) \cong \pi_1(C) \times \pi_1(C) \lra 1
\]
where $R$ is the normal closure in $\Gam$ of the torsion elements associated with the orbifold locus.

Genus two curves admit \'etale covers by a curve of any given genus greater than two, so given $g_1, g_2 \ge 2$ there exist curves $C_j$ of genus $g_j$ and \'etale covers $C_j \to C$. In particular, there is an \'etale cover $C_1 \times C_2 \to C \times C$. This determines a finite index subgroup $\Del^\prime < \pi_1(Y)$ and a cover $Y^\prime \to Y$ with $Y^\prime$ birational to $C_1 \times C_2$. The pullback of $\Del^\prime$ to $\Gam$ under the surjection is a finite index subgroup $\Gam^\prime$ with the property that $\Gam^\prime \bs \bbB^2$ is an orbifold with underlying analytic space $Y^\prime$. This proves the proposition.
\end{pf}

\begin{rem}\label{rem:Higherg}
The surface birational to $C_1 \times C_2$ that is biholomorphic to a ball quotient in the proof of Proposition \ref{prop:MainHyp} is the blowup of $C_1 \times C_2$ at $3(g_1-1)(g_2-1)$ points.
\end{rem}

The next proposition handles higher-genus ruled surfaces in a similar fashion to Proposition \ref{prop:MainRuled}, again building off of the example from \S \ref{sec:Ex1}.

\begin{prop}\label{prop:MainRuled}
Given $g_2 \ge 2$, there exists a smooth projective curve $C_2$ of genus $g_2$ so that $\bbP^1 \times C_2$ is birational to a ball quotient.
\end{prop}

\begin{pf}
Let $C \times C$ be as in Theorem \ref{thm:Ex}, and consider the action of $\bbZ / 5$ on the first coordinate by $\xi$. The quotient $C / \langle \xi \rangle$ is isomorphic to $\bbP^1$. Since $\xi$ fixes the three points $a_j \in \{\pm 1, \infty\}$, it fixes the points $x_j = (a_j, a_j)$ on the product. This implies that there is an induced action of $\xi$ on the blowup surface $Y$. Moreover, $\xi$ cyclically permutes the graphs $C_j$ of $\xi^j$, and thus the action on $Y$ permutes the curves $D_j$ that form the orbifold locus of $\Gam \bs \bbB^2$ for the ball quotient structure on $Y$.

If $\Aut_\calD(Y)$ denotes automorphisms of $Y$ preserving $\calD$, recall from the discussion around Lemma \ref{lem:Gquotient} that the quotient $Y / \Aut_\calD(Y)$ is again a ball quotient. This implies that the action of $\xi$ induces an action of $\bbZ / 5$ on $\Gam \bs \bbB^2$, meaning that there is a lattice $\Gam_\xi < \PU(2,1)$ containing $\Gam$ as a normal subgroup of index $5$ so that the orbifold covering $\Gam \bs \bbB^2 \to \Gam_\xi \bs \bbB^2$ is induced by this action.

Now, the complement $Y^\circ = Y \ssm \bigcup E_j$ of the exceptional curves in $Y$ is isomorphic to $C \times C$ minus the three points $x_j$. Thus $Y^\circ / (\bbZ / 5)$ is isomorphic to $\bbP^1 \times C = (C \times C) / \langle \xi \rangle$ minus the three points $(\infty, a_0)$, $(1, a_1)$, and $(-1, a_2)$. In other words, the underlying orbifold for $\Gam_\xi \bs \bbB^2$ is birational to $\bbP^1 \times C$. This proves the proposition for the case $g_2 = 2$. The case $g_2 > 2$ is now proved using \'etale covers of $C$ exactly as in Proposition \ref{prop:MainHyp}.
\end{pf}

\begin{rem}
The surfaces in the proof of Proposition \ref{prop:MainRuled} biholomorphic to a ball quotient are singular surfaces birational to a smooth ruled surface. On the other hand, it follows from the proof that for all $g \ge 2$ there is a curve $C_g$ of genus $g$ such that the smooth surface $\bbP^1 \times C_g$ is holomorphically dominated by a ball quotient.
\end{rem}

The following lemma will be used to construct examples for the case $g_1 = 1$ and $g_2 \ge 2$.

\begin{lem}\label{lem:TriangleSubgroup}
The curve $C$ from \S \ref{sec:Ex1} has an \'etale cover $C^\prime$ that admits a holomorphic mapping onto an elliptic curve $E^\prime$.
\end{lem}

\begin{pf}
The fact that $C$ is the cyclic cover of $\bbP^1$ branched over three points to order five means that the uniformization of $C$ as a hyperbolic $2$-manifold $\Sig \bs \bbH^2$ realizes $\Sig$ as an index five subgroup of the triangle group $\Del(5,5,5)$, which itself is an index six subgroup of $\Del(2,3,10)$. From the presentation
\[
\Del(2,3,10) =\! \left\langle a,b\ |\ a^2, b^3, (ab)^{10}\right\rangle
\]
one can show that the subgroup $\Lam$ generated by the words $(ab)^2 a b^{-1} aba$ and $abab^{-1} a (ba)^2$ has index ten with quotient $\Lam \bs \bbH^2$ having genus zero and one cone point of order three. Thus the quotient is isomorphic to an elliptic curve $E^\prime$. Letting $C^\prime$ be the curve with fundamental group $\Sig \cap \Lam$ that completes the diagram
\[
\begin{tikzcd}[column sep = small]
(\Sig \cap \Lam) \bs \bbH^2 = C^\prime \arrow[dd] \arrow[drr] & & \\
& & \Sig \bs \bbH^2 = C \arrow[dd] \\
\Lam \bs \bbH^2 \cong E^\prime \arrow[ddr] & & \\
 & & \Del(5,5,5) \bs \bbH^2 \cong \bbP^1 \arrow[dl] \\
 & \Del(2,3,10) \bs \bbH^2 \cong \bbP^1 &
\end{tikzcd}
\]
gives the desired cover.
\end{pf}

\begin{prop}\label{prop:Elliptic}
Let $C$ be the curve with genus $2$ from \S \ref{sec:Ex1} and $C^\prime \to E^\prime$ be the map from an \'etale cover of $C$ to an elliptic curve $E^\prime$ provided by Lemma \ref{lem:TriangleSubgroup}. For all $g \ge 2$ there is an \'etale cover $C_g$ of $C$ and a compact ball quotient $\Gam_g \bs \bbB^2$ that admits a dominant rational map onto $E^\prime \times C_g$.
\end{prop}

\begin{pf}
From the proof of Proposition \ref{prop:MainHyp}, there is a ball quotient $\Gam_g \bs \bbB^2$ with underlying analytic space isomorphic to the blowup of $C^\prime \times C_g$ at some collection of points, where $C_g$ is any \'etale cover of $C$ with genus $g$. Contracting the exceptional curves an applying the mapping $C^\prime \to E^\prime$ on the first coordinate proves the proposition.
\end{pf}

The remaining example required to prove Theorem \ref{thm:Main} is constructed using the example from \S \ref{sec:Ex2}.

\begin{prop}\label{prop:P1E}
There is a compact ball quotient $\Gam^\prime \bs \bbB^2$ whose underlying analytic space is biholomorphic to a ruled surface of genus one.
\end{prop}

\begin{pf}
The proof will employ the example and notation from \S \ref{sec:Arith2}. In particular, let $\Gam \bs \bbB^2$ be the example from \S \ref{sec:Ex2} and consider the automorphism
\[
\beta = \begin{pmatrix} 1 & 0 \\ 1 & \zeta^2 \end{pmatrix}
\]
of $T \times T$ with order $3$. This acts on $T \times T$ as a complex reflection stabilizing the curve $T_\beta$ parametrized by
\[
z \longmapsto\!\left((1-\zeta^2) z\,,\, z\right)\!.
\]
The induced action on the blowup $Y$ is again by a complex reflection of order $3$ (the multiplicative order of $\zeta^2$) through the proper transform $D_\beta$ of $T_\beta$ to $Y$.

Note that $\beta$ acts trivially on the first factor of $T \times T$ and for a fixed $z_0$ acts on the vertical curve $(z_0, w)$ by $w \mapsto \zeta(z_0 - w)$, which has quotient $\bbP^1$. Thus $(T \times T) / \langle \beta \rangle$ is simply $T \times \bbP^1$, and this realizes $T \times T$ as the $3$-fold branched cover of $T \times \bbP^1$ branched over the graph of multiplication by $\zeta^2$. It was shown in \S \ref{sec:Arith2} that $\beta$ acts on $\Gam \bs \bbB^2$, and it follows exactly as in the proof of Proposition \ref{prop:MainRuled} that $(\Gam \bs \bbB^2) / \langle \beta \rangle$ is a ball quotient with underlying analytic space birational to $\bbP^1 \times T$. This proves the proposition.
\end{pf}

The results in this section now prove Theorem \ref{thm:Main}.

\begin{pf}[Proof of Theorem \ref{thm:Main}]
Combine Propositions \ref{prop:MainHyp}, \ref{prop:MainRuled}, \ref{prop:Elliptic}, and \ref{prop:P1E}.
\end{pf}

\section{Surfaces with $\kappa(V) = 0$}\label{sec:Kodaira}

The goal of this section is to prove Theorem \ref{thm:Kodaira}. By Corollary \ref{cor:MainCor}, it suffices to consider surfaces of Kodaira dimension $0$. The case of an abelian surface is covered by the example in \S \ref{sec:Ex2}. By \cite[Thm.\ VIII.2]{Beauville}, the remaining smooth projective surfaces of Kodaira dimension zero are:
\begin{itemize}
\item[$\star$] K3 surfaces, which define one connected deformation space \cite[Thm.\ 22.3]{BPV}
\item[$\star$] Enriques surfaces, which form one connected deformation space \cite[Thm.\ 18.5]{BPV}
\item[$\star$] bielliptic surfaces, which fall into seven connected deformation spaces \cite[List VI.20]{Beauville}, \cite[Rem.\ 2]{Suwa}
\end{itemize}
Bielliptic surfaces are surfaces of the form $(E_\lam \times E_\tau) / G$, where $E_\lam$ and $E_\tau$ are elliptic curves and $G$ is a finite group acting freely on the product. By the Bagnera--de Franchis theorem \cite[List VI.20]{Beauville},
\[
G \in\!\left\{\bbZ / 2, (\bbZ / 2)^2, \bbZ / 4, \bbZ / 4 \oplus \bbZ / 2, \bbZ / 3, (\bbZ / 3)^2, \bbZ / 6 \right\}\!,
\]
where each group determines exactly one connected deformation space of surfaces \cite[Rem.\ 2]{Suwa}. As described in more detail below, each deformation space of surfaces of Kodaira dimension zero contains an element somehow closely related to a product of elliptic curves, and this connection will be used to prove Theorem \ref{thm:Kodaira}.

A second abelian example is required to cover all bielliptic surfaces. This is a compact variant on a noncompact example due to Holzapfel \cite{HolzapfelAbelian} also considered in \cite{DiCerboStoverGauss}. The example presented here is analogous to the example in \S\ref{sec:Ex2} alongside Hirzebruch's noncompact ball quotient from the same paper (which appears in the proof of Proposition \ref{prop:FoundDM}) based on the abelian surface defined using the Eisenstein integers. Consider the elliptic curve $E = \bbC / \bbZ[i]$ defined by the Gaussian integers and the curves
\begin{align*}
C_0 &= (z,z) & H_0 &=\!\left(\!z\,,\, \frac{1}{2}\!\right) \\
C_1 &= (z,-z) & H_1 &=\!\left(\!z\,,\, \frac{i}{2}\!\right) \\
C_2 &= (z,iz) & V_0 &=\!\left(\!\frac{1}{2}\,,\, z\!\right) \\
C_3 &= (z,-iz) & V_1 &=\!\left(\!\frac{i}{2}\,,\, z\!\right)
\end{align*}
on $E \times E$. These eight curves intersect in six points, with exactly four of the curves passing through each point and exactly four intersection points on each $C_j$ and two on each $H_j$ and $V_j$. See \cite[\S 3]{DiCerboStoverGauss} for the coordinates of the intersection points.
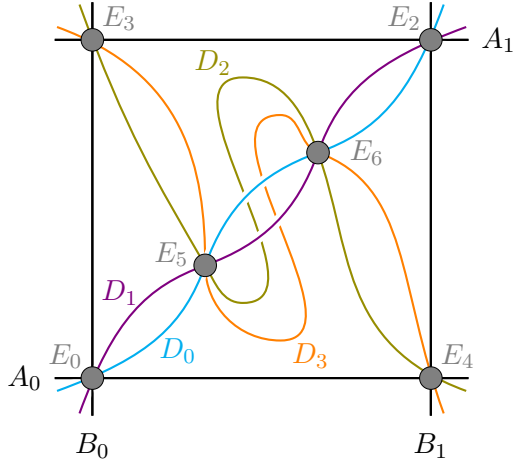
\begin{figure}[h]
\centering
\begin{tikzpicture}[scale=1.5]
\draw[thick, shorten >= -0.5cm, shorten <= -0.5cm] (0,0) -- (0,3);
\draw[thick, shorten >= -0.5cm, shorten <= -0.5cm] (0,3) -- (3,3);
\draw[thick, shorten >= -0.5cm, shorten <= -0.5cm] (3,3) -- (3,0);
\draw[thick, shorten >= -0.5cm, shorten <= -0.5cm] (3,0) -- (0,0);
\draw[thick, olive, shorten >= -0.5cm, shorten <= -0.5cm] (0,3) to[out=290, in=120] (1,1) to[out=300, in=180] (1.3333, 0.6666) to[out=0, in=180] (1.3333, 2.6666) node[xshift=-0.4cm, yshift=0.2cm] {$D_2$} to[out=0, in=110] (2,2) to[out=290, in=160] (3,0);
\draw[thick, orange, shorten >= -0.5cm, shorten <= -0.5cm] (0,3) to[out=-20, in=90] (1,1) to[out=270, in=180] (1.6666, 0.3333) node[xshift=0.4cm, yshift=-0.2cm] {$D_3$} to[out=0, in=180] (1.6666, 2.3333) to[out=0, in=160] (2,2) to[out=-20, in=110] (3,0);
\draw[line width=4pt, white, shorten >= -0.5cm, shorten <= -0.5cm] (0,0) to[out=70, in=200] (1,1) to[out=20, in=250] (2,2) to[out=70, in=200] (3,3);
\draw[line width=4pt, white, shorten >= -0.5cm, shorten <= -0.5cm] (0,0) to[out=20, in=250] (1,1) to[out=70, in=200] (2,2) to[out=20, in=250] (3,3);
\draw[thick, violet, shorten >= -0.5cm, shorten <= -0.5cm] (0,0) to[out=70, in=200] node[xshift=-0.2cm, yshift=0.2cm] {$D_1$} (1,1) to[out=20, in=250] (2,2) to[out=70, in=200] (3,3);
\draw[thick, cyan, shorten >= -0.5cm, shorten <= -0.5cm] (0,0) to[out=20, in=250] node[xshift=0.2cm, yshift=-0.2cm] {$D_0$} (1,1) to[out=70, in=200] (2,2) to[out=20, in=250] (3,3);
\draw[thick, shorten >= -0.5cm, shorten <= -0.5cm] (0,0) node[xshift=-0.9cm] {$A_0$} node[yshift=-0.9cm] {$B_0$} -- (0,3);
\draw[thick, shorten >= -0.5cm, shorten <= -0.5cm] (0,3) -- (3,3);
\draw[thick, shorten >= -0.5cm, shorten <= -0.5cm] (3,3) node[xshift=0.9cm] {$A_1$} -- (3,0);
\draw[thick, shorten >= -0.5cm, shorten <= -0.5cm] (3,0) node[yshift=-0.9cm] {$B_1$} -- (0,0);
\draw[very thin, fill=gray] (0,0) circle (0.1cm) node[above left, gray] {$E_0$};
\draw[very thin, fill=gray] (3,0) circle (0.1cm) node[above right, gray] {$E_4$};
\draw[very thin, fill=gray] (0,3) circle (0.1cm) node[above right, gray] {$E_3$};
\draw[very thin, fill=gray] (3,3) circle (0.1cm) node[above left, gray] {$E_2$};
\draw[very thin, fill=gray] (1,1) circle (0.1cm) node[xshift=-0.45cm, yshift=0.15cm, gray] {$E_5$};
\draw[very thin, fill=gray] (2,2) circle (0.1cm) node[xshift=0.65cm, gray] {$E_6$};
\end{tikzpicture}
\caption{The blown-up Gaussian abelian surface}\label{fig:Gauss}
\end{figure}

Let $Y$ be the blowup of $E \times E$ at these six points, $D_j$ (resp.\ $A_j$, $B_j$) be the proper transform of $C_j$ (resp.\ $H_j$, $V_j$) and $E_j$ be the exceptional curves, ordered as indicated in Figure \ref{fig:Gauss}. Then consider the $\bbQ$-divisor
\begin{equation}\label{eq:GaussD}
\calD =\! \left(1 - \frac{1}{2}\right) \!\sum E_j +\! \left(1 - \frac{1}{4}\right)\!\left(\sum A_j + \sum B_j + \sum D_j \right)
\end{equation}
on $Y$. Note that the curves $A_j$, $B_j$, and $D_j$ are all disjoint, as are the six exceptional curves, and each $E_j$ meets exactly four curves. Then
\begin{alignat*}{2}
c_1^2(Y, \calD) &= 6(-1)\!\left(2 - \frac{1}{2}\right)^2\! &&+ 4(-2)\!\left(1 - \frac{1}{4}\right)^2\! + 4(-4)\!\left(1 - \frac{1}{4}\right)^2\! \\
& &&+ 2(6)(4) \!\left(2 - \frac{1}{2}\right)\!\left(1 - \frac{1}{4}\right) \\
&= 27 &&
\end{alignat*}
since $K_Y$ is numerically equivalent to $\sum E_j$, $E_j^2 = -1$, $A_j^2 = B_j^2 = -2$, and $D_j^2 = -4$. As for $c_2(Y, \calD)$, $c_2(Y) = 6$ and every normal crossing of two curves has local group $\bbZ / 4 \times \bbZ / 2$, thus weight $8$, and so
\begin{alignat*}{2}
c_2(Y, \calD) &= 6 &&- 6(-2)\!\left(1 - \frac{1}{2}\right)\! - 4(-2)\!\left(1 - \frac{1}{4}\right) \\
& &&- 4(-4)\!\left(1 - \frac{1}{4}\right)\! - 24\! \left(1 - \frac{1}{8}\right) \\
&= 9
\end{alignat*}
and hence $c_1^2(Y, \calD) = 3 c_2(Y, \calD)$. The proof that $K_Y + \calD$ is ample is exactly the same as in \S\ref{sec:Ex2}: one only needs to check that $K_Y + \calD$ pairs positively with the curves in the arrangement under the intersection pairing. This proves the following result.

\begin{prop}\label{prop:GaussEx}
Let $Y$ be the blown up abelian surface defined in this section and $\calD$ the $\bbQ$-divisor on $Y$ defined in Equation \eqref{eq:GaussD}. Then $(Y, \calD)$ is a ball quotient pair.
\end{prop}

\begin{rem}\label{rem:NA1}
Giving each curve in the support of $\calD$ orbifold weight $3$ also defines a ball quotient pair. This alternate example is actually nonarithmetic; see Remark \ref{rem:NA2} below.
\end{rem}

The example from Proposition \ref{prop:GaussEx} and the example from \S \ref{sec:Ex2} will now be used to cover to other surfaces of Kodaira dimension zero and prove Theorem \ref{thm:Kodaira}.

\begin{pf}[Proof of Theorem \ref{thm:Kodaira}]
The proof proceeds by constructing an example from each deformation space. Note that the abelian surface case is handled by \S \ref{sec:Ex2} and Proposition \ref{prop:GaussEx}.

\bigskip

First, let $(Y, \calD)$ be the ball quotient pair from Proposition \ref{prop:GaussEx} and $\Gam \bs \bbB^2$ the associated ball quotient orbifold. This will be used to produce a K3 surface, an Enriques surface, and all but three types of bielliptic surfaces.

\bigskip

\noindent
\emph{K3 and Enriques surfaces.} The union of the curves $C_j$, $H_j$, and $V_j$ on $E \times E$ used to define $\calD$ is stable under multiplication by $i \Id$ in such a way that orbifold weights are preserved, which implies that it induces an action of $\bbZ / 4$ on $\Gam \bs \bbB^2$. If $\Del$ is the associated lattice in $\PU(2,1)$ containing $\Gam$ as an index four subgroup, then $\Del \bs \bbB^2$ is birational to $(E \times E) / \langle i \Id \rangle$. Let $\Lam < \Del$ be the intermediate subgroup associated with $\langle - \Id \rangle$. Then $(E \times E) / \langle - \Id \rangle$ is a Kummer variety \cite[V.16]{BPV}, and thus $\Lam \bs \bbB^2$ is birational to a K3 surface. Moreover, the induced action of $\bbZ / 2$ on $\Lam \bs \bbB^2$ is free away from the singular points of the Kummer surface, hence $\Del \bs \bbB^2$ is birational to an Enriques surface \cite[Prop.\ VIII.17]{Beauville}.

\bigskip

The bielliptic examples constructed from the pair $(Y, \calD)$ from Proposition \ref{prop:GaussEx} are all compact analogues of noncompact examples in \cite[\S3]{DiCerboStoverGauss}.

\bigskip

\noindent
\emph{A $\bbZ / 2$ bielliptic surface.} Consider the automorphism
\[
\psi_2(z, w) =\! \left(z + \frac{1+i}{2}\,,\,-w + \frac{1+i}{2} \right)
\]
of $E \times E$ with order two. This automorphism permutes the curves $C_j$, $H_j$, and $V_j$ and hence permutes their intersection points and induces an automorphism $\wh{\psi}_2$ of $Y$. Moreover, the action on $Y$ fixes the collection of curves with the same orbifold weight on $\Gam \bs \bbB^2$, thus there is a lattice $\Lam_2$ containing $\Gam$ with index two so that $\Lam_2 \bs \bbB^2$ has underlying analytic space $Y / \langle \wh{\psi}_2\rangle$. Conjugating $\psi_2$ by translation by $(1+i)/4$ in the second coordinate shows that it is equivalent to the automorphism
\[
\phi_2(z, w) =\! \left(z + \frac{1+i}{2}\,,\,-w \right),
\]
which is the standard automorphism defining a $\bbZ / 2$ bielliptic quotient of $E \times E$. Thus $\Lam_2 \bs \bbB^2$ is birational to a $\bbZ / 2$ bielliptic surface.

\bigskip

\noindent
\emph{A $\bbZ / 4$ bielliptic surface.} Consider the $4$-fold \'etale cover $E^\prime \times E$ of $E \times E$ associated with the subgroup of $\bbZ[i]$ generated by $4$ and $i$. This induces a ball quotient orbifold $\Gam^\prime \bs \bbB^2$ covering $\Gam \bs \bbB^2$ with underlying analytic space the blowup $Y^\prime$ of $E^\prime \times E$ at $24$ points and orbifold locus the preimage $\calD^\prime$ of $\calD$. The preimage $C_j^\prime$ of each curve $C_j$ in $E^\prime \times E$ is irreducible, and can be described in coordinates in exactly the same way by considering $C_j^\prime$ as the image in $E^\prime \times E$ of the graph of multiplication by $i^j$ in $\bbC^2$. The curves $H_j$ similarly have irreducible preimage $H_j^\prime$ given as the image of the appropriate horizontal $\bbC$ in $\bbC^2$, but the curves $V_j$ each lift to four irreducible curves whose coordinates are of the form $(\tau, z)$ where $\tau$ is one of the four preimages of $1/2, i/2 \in E$ in $E^\prime$. Thus the support of $\calD^\prime$ consists of fourteen curves.

The point $1 \in E^\prime$ is a $4$-torsion point that generates the group of deck transformations for the covering map to $E$. Therefore
\[
\psi_4(z, w) = (z + 1, i w)
\]
is an automorphism of $E^\prime \times E$ with quotient a $\bbZ / 4$ bielliptic surface. Direct check shows that $\psi_4$ preserves $\calD^\prime$ and hence induces an automorphism of $\Gam^\prime \bs \bbB^2$ with quotient birational to a $\bbZ / 4$ bielliptic surface.

\bigskip

\noindent
\emph{A $(\bbZ / 2)^2$ and a $\bbZ / 4 \times \bbZ / 2$ bielliptic surface.} Let $\Gam_4 \bs \bbB^2$ be the ball quotient in this section that is birational to a $\bbZ / 4$ bielliptic surface. As shown at the end of the proof of \cite[Thm.\ 1.4]{DiCerboStoverGauss}, there are \'etale covers of a $\bbZ / 4$ bielliptic surface by both $(\bbZ / 2)^2$ and $\bbZ / 4 \times \bbZ / 2$ bielliptic surfaces. This induces orbifold covers of $\Gam_4 \bs \bbB^2$ with underlying analytic space birational to surfaces of the desired kind, hence completing these two cases.

\bigskip

The remainder of the proof relies on the example from \S\ref{sec:Ex2}, which was constructed using the Eisenstein integers $\bbZ[\zeta]$, where $\zeta$ is a primitive $6^{th}$ root of unity. Recall that $T$ is the elliptic curve $\bbZ[\zeta] \bs \bbC$. The arguments are completely analogous to the previous ones, except the surfaces are now closely related to the bielliptic ball quotients in \cite[\S 6]{DiCerboStoverClassify} instead of those found in \cite{DiCerboStoverGauss}. While these examples are commensurable with the example in \S\ref{sec:Ex2}, it is easier to describe them explicitly instead of giving various group actions in coordinates. It will be convenient to let $\rho = \zeta^2$ be a primitive $3^{rd}$ root of unity, and recall that $\rho-1$ divides $3$.

\bigskip

\noindent
\emph{A $\bbZ / 3$ bielliptic surface.} This example is analogous to \cite[\S 6.3]{DiCerboStoverClassify}. The automorphism
\[
\phi_3(z, w) =\! \left(z + \frac{\rho-1}{3}\,,\, \rho w\right)
\]
of the abelian surface $T \times T$ from \S\ref{sec:Ex2} has order three and quotient a $\bbZ / 3$ bielliptic surface $Z_3$. First, notice that the image $C_0^\prime$ in $Z_3$ of the diagonal $C_0$ in $T \times T$ has genus one and a unique singular point of order three. Indeed, $\phi_3(C_0) \cap C_0$ consists of the $\phi_3$-orbit of the $3$-torsion point $(1/3, 1/3)$, which therefore maps to a singular point of $C_0^\prime$. If $C_1$ is the horizontal curve $(z, 1/3)$ on $T \times T$, then $\phi_3(C_1)$ does not meet $C_1$, hence the image $C_1^\prime$ is a smooth curve of genus one passing through the singular point of $C_0^\prime$.

Let $Y_3$ be the blowup of $Z_3$ at the singular point of $C_0^\prime$ with exceptional curve $E$ and $D_j$ be the proper transform of the image of $C_j^\prime$, which is a smooth curve on $Y_3$. Then $D_0^2 = -1$ and $D_1^2 = -3$ by adjunction \cite[I.15]{Beauville}, since $E$ is numerically equivalent to the canonical divisor, $D_0 \cdot E = 1$, and $D_1 \cdot E = 3$. Then
\[
\calD_3 =\!\left(1- \frac{1}{3}\right)\!(E + D_0 + D_1)
\]
has $c_1^2(Y_3, \calD_3) = 13/3$ and $c_2(Y_3, \calD_3) = 13/9$. A direct check analogous to all other cases in this paper shows that $K_{Y_3} + \calD$ is ample, and thus $(Y_3, \calD_3)$ is a ball quotient pair by Theorem \ref{thm:OrbiBMY}.

\bigskip

\noindent
\emph{A $(\bbZ / 3)^2$ bielliptic surface.} This example is analogous to \cite[\S 6.4]{DiCerboStoverClassify}. Consider the $3$-fold cover $T^\prime$ of $T$ associated with the index $3$ subgroup $(\rho - 1) \bbZ[\zeta]$ of $\bbZ[\zeta]$. Then the $3$-torsion subgroup of $T^\prime$ is generated by $1$ and $(\rho - 1)/3$. The automorphisms
\begin{align*}
\psi_a(z, w) &= (z + 1, \rho w) & \psi_b(z, w) =\! \left(z + \frac{\rho - 1}{3}, w + \frac{\rho - 1}{3}\right)
\end{align*}
of $T^\prime \times T$ commute, and the quotient by $T^\prime \times T$ by the group they generate is a $(\bbZ / 3)^2$ bielliptic surface $Z_9$. As in the previous example, if $C_0$ is the image of the diagonal in $\bbC^2$ in $T^\prime \times T$ and $C_0^\prime$ the image of $C_0$ in $Z_9$, then $C_0^\prime$ has genus zero and one singular point of order $3$. Again, let $C_1^\prime$ be the horizontal curve on $Z_9$ through the singular point of $C_0^\prime$.

Let $Y_9$ be the blowup of $Z_9$ at the singular point of $C_0^\prime$ with exceptional curve $E$ and $D_j$ be the proper transform of $C_j^\prime$ to $Y_9$. Then the triple $E, D_0, D_1$ has exactly the same intersection properties as the curves in the previous $\bbZ / 3$ bielliptic example. It therefore follows by the exact same calculations that if
\[
\calD_9 =\!\left(1- \frac{1}{3}\right)\!(E + D_0 + D_1)
\]
then $(Y_9, \calD_9)$ is a ball quotient pair.

\bigskip

\noindent
\emph{A $\bbZ / 6$ bielliptic surface.} To complete the proof of the theorem, it remains to produce a lattice $\Gam < \PU(2,1)$ so that $\Gam \bs \bbB^2$ admits a dominant rational map to a $\bbZ / 6$ bielliptic surface. Let $\Gam$ be a lattice so that the underlying analytic space $Y_3$ is the blowup of a $\bbZ / 3$ bielliptic surface $Z_3$ at a single point. Any $\bbZ / 6$ bielliptic surface is the quotient of a $\bbZ / 3$ bielliptic by a free action of $\bbZ / 2$, and every $\bbZ / 3$ bielliptic admits such a free action. In particular, there is a $2$-to-$1$ \'etale map from $Z_3$ to some bielliptic surface $Z_6$. The composition
\[
\Gam \bs \bbB^2 \lra Y_3 \lra Z_3 \lra Z_6
\]
is dominant and rational. This completes the proof of the theorem.
\end{pf}

\begin{rem}\label{rem:Kasparian}
Holzapfel noticed that the complement in $Y$ of the curves $A_j, B_j, D_j$ is a noncompact ball quotient \cite{HolzapfelAbelian}, and that one could take a quotient by the action of $-\Id$ to get a ball quotient with compactification birational to a K3 surface. Kasparian and Kotzev noticed that the further quotient by $i \Id$ gives a ball quotient with compactification birational to an Enriques surface \cite[Thm.\ 9]{KasparianKotzev}. The K3 and Enriques cases in the proof of Theorem \ref{thm:Kodaira} are nothing more than applying the exact same argument to the compact ball quotient orbifold from Proposition \ref{prop:GaussEx}.
\end{rem}

\begin{rem}
Careful analysis shows that the $\bbZ / 6$ bielliptic surface in the proof of Theorem \ref{thm:Kodaira} does not inherit a ball quotient orbifold structure, and this is the only case in Kodaira dimension zero where the proof does not provide a ball quotient orbifold structure on a surface birational to the desired smooth minimal model. It would be interesting to settle this case like the others by an explicit ball quotient orbifold structure.
\end{rem}

\begin{rem}
If one uses toroidal compactifications instead of compact ball quotients, almost every case in Theorem \ref{thm:Kodaira} was in the literature already. The case of an abelian surface is covered by the examples in \S \ref{sec:Ex2}, either by the compact example or the auxiliary example (also due to Hirzebruch \cite{HirzebruchChern}) that admits a smooth toroidal compactification by the four elliptic curves $D_j$. See \cite{DiCerboStoverMultiple} for a plethora of other related (albeit commensurable) examples. See Remark \ref{rem:Kasparian} for the K3 and Enriques cases. Arithmetic lattices in $\PU(2,1)$ with smooth toroidal compactification a bielliptic surface associated with $G$ either $\bbZ / 3$ or $(\bbZ / 3)^2$ were constructed in \cite{DiCerboStoverClassify, DiCerboStoverComm} and examples for the groups $\bbZ / 2$, $(\bbZ / 2)^2$, $\bbZ / 4$, $\bbZ / 4 \oplus \bbZ / 2$ were constructed in \cite{DiCerboStoverGauss}. In other words, the only surface of Kodaira dimension zero not already constructed in the literature in the more general setting of toroidal compactifications is a bielliptic surface associated with the group $\bbZ / 6$.
\end{rem}

The final task in this section is to prove Theorem \ref{thm:Farb}.

\begin{pf}[Proof of Theorem \ref{thm:Farb}]
The main example in \S \ref{sec:Ex2} is a ball quotient orbifold with underlying space the blowup of an abelian surface at a single point. This smoothly dominates any abelian surface. Similarly, the proof of Theorem \ref{thm:Kodaira} constructs ball quotient orbifolds with underlying space the blowup of a bielliptic surface, and these examples suffice to show that any bielliptic surface is smoothly dominated by a ball quotient. Thus only K3 and Enriques surfaces remain.

Naruki \cite{Naruki} constructed a K3 ball quotient as follows. See Ch.\ VIII and \S V.7 of \cite{BPV} for basic information about the structure of K3 surfaces and elliptic fibrations used in what follows. There exists a K3 surface admitting an elliptic fibration $X \to \bbP^1$ with the following properties:
\begin{itemize}

\item[$\star$] There are exactly three singular fibers of type $I_7$ and exactly three of type $I_1$.

\item[$\star$] The fibration admits seven distinct disjoint sections for which each section intersects exactly one irreducible component of each $I_7$ singular fibers and each curve on each $I_7$ singular fiber meets exactly one section.

\end{itemize}
Let $A_{j, k}$ denote the $k^{th}$ curve in the $j^{th}$ singular fiber of type $I_7$ and $B_k$ be the $k^{th}$ section (for any ordering). Then $A_{j,k}^2 = B_k^2 = -2$, and any curve meets exactly three other curves.

Consider the $\bbQ$-divisor
\[
\calC =\! \left(1 - \frac{1}{7}\right)\! \sum_{k=1}^7\!\left(B_k + \sum_{j = 1}^3 A_{j,k} \right)
\]
on $X$. Since K3 surfaces have trivial canonical bundle and $c_2(X) = 24$,
\begin{align*}
c_1^2(X, \calC) &=\!\left(\frac{6}{7}\right)^2\!\left(28(-2) + 2\!\left(\frac{28}{2}(3)\right)\!\right) \\
&= \frac{144}{7} \\
c_2(X, \calC) &= 24 - (28)\!\left(\frac{6}{7}\right)\!(-1) - (42)\!\left(\frac{48}{49}\right)\! \\
&= \frac{48}{7}
\end{align*}
hence $c_1^2(X, \calC) = 3 c_2(X, \calC)$. Moreover, $\calC$ is ample. Indeed, any reduced irreducible curve not in the support of $\calC$ must meet either the $I_7$ singular fiber or some section. Then each curve in the support of $\calC$ meets $\calC$ with intersection number $(6/7)(-2 + 3) = 6/7$, so $\calC$ is ample by the Nakai--Moishezon criterion \cite[Thm.\ V.1.10]{Hartshorne}. Thus $(X, \calC)$ is a ball quotient pair, and it is arithmetic by \cite[Thm.\ 1]{Naruki}.

Since all K3 surfaces are diffeomorphic \cite[Cor.\ VII.8.6]{BPV}, all K3 surfaces are dominated by $X$, and therefore all K3 surfaces are smoothly dominated by a ball quotient. Since every Enriques surface has universal cover a K3 surface, all Enriques surfaces are thus also smoothly dominated by a ball quotient. This completes the proof of the theorem.
\end{pf}

\section{Arithmeticity of the examples in \S \ref{sec:Kodaira}}\label{sec:LastArithmetic}

To prove Theorem \ref{thm:Arithmetic}, what remains is to prove arithmeticity of the examples from \S \ref{sec:Kodaira} that are not arithmetic by construction.

\begin{thm}\label{thm:GaussExA}
The ball quotient pair $(Y, \calD)$ from Proposition \ref{prop:GaussEx} is arithmetic.
\end{thm}

\begin{pf}
Recall that $E$ is the elliptic curve associated with the Gaussian integers. Consider the action of $F = (\bbZ / 4)^2$ on $E \times E$ generated by complex multiplication by $i$ on each factor. First, note that $F$ stabilizes the set of points blown up to obtain $Y$. Indeed, the six points are
\begin{align*}
(0,0) && \left(\frac{1}{2}, \frac{1}{2}\right) && \left(\frac{i}{2}, \frac{i}{2}\right) && \left(\frac{1+i}{2}, \frac{1+i}{2}\right) \\
&& \left(\frac{1}{2}, \frac{i}{2}\right) && \left(\frac{i}{2}, \frac{1}{2}\right) && 
\end{align*}
where $(0,0)$ and $((1+i)/2, (1+i)/2)$ are fixed by all of $F$ and the other four are permuted with stabilizer $(\bbZ / 2)^2$.

The remaining points of $E \times E$ fixed by some nontrivial element of $F$ are of the form $(\tau_1, \tau_2)$ for $\tau_j$ a $2$-torsion point, and each stabilizer is either $\bbZ / 4 \times \bbZ / 2$ or $\bbZ / 2 \times \bbZ / 4$. Finally, note that $F$ permutes the curves $C_j$, each having stabilizer $\bbZ / 4$ acting as complex multiplication. In particular, the induced action on $Y$ permutes curves in the support of $\calD$ with the same orbifold weight, hence $Z = Y / F$ is a ball quotient orbifold with respect to the appropriate $\bbQ$-divisor $\calE$ with support the image in $Z$ of the orbifold locus for $(Y, \calD)$ and any additional curves fixed by the action of $F$.

Since $(E \times E) / F$ is $\bbP^1 \times \bbP^1$ and the six exceptional curves of $Y$ map to three exceptional curves of $Z$, it follows that $Z$ is the blowup of $\bbP^1 \times \bbP^1$ on three points on the diagonal. As mentioned in the proof of Lemma \ref{lem:ZisP2}, this means that $Z$ is isomorphic to the blowup of $\bbP^2$ at the vertices of the complete quadrangle. Moreover, the curves in $Y$ stabilized by some subgroup of $F$ are mapped to the lines in the quadrangle.

It remains to compute orbifold weights of the fixed curves and compare the orbifold structure on $Z$ with Deligne--Mostow orbifolds. First, the four curves $D_j$ are identified by $F$, and the action of the stabilizer is by complex multiplication of order four. In particular, the image of all the $D_j$ in $Z$ is a single curve of orbifold weight four. The points $(0,0)$ and $((1+i)/2, (1+i)/2)$ have stabilizer $F$, where the diagonal subgroup of $F$ therefore acts as a complex reflection of order four through the associated exceptional curves, which are $E_5$ and $E_6$ in Figure \ref{fig:Gauss}. Since each $E_j$ has orbifold weight two, the image on $Z$ has orbifold weight eight. The remaining blowup points $(1/2, 1/2)$, $(1/2, i/2)$, $(i/2, 1/2)$, $(i/2, i/2)$ are permuted by $F$ and have stabilizer $(\bbZ/2)^2$, so the diagonal subgroup acts as a complex reflection of order two through $E_j$ for $1 \le j \le 4$, and hence the image in $Z$ is a single curve of orbifold weight four.

The curves $A_j$ (resp.\ $B_j$) coming from horizontal and vertical curves are exchanged by one direct $\bbZ / 4$ factor of $F$ with $\bbZ / 2$ stabilizer acting as a complex reflection, and the other direct $\bbZ / 4$ factor of $F$ acts by complex multiplication. Since the $A_j$ and $B_j$ have orbifold weight four, their images on $Z$ determine two curves in total, each with orbifold weight eight. Similarly, the horizontal curves $(z, 0)$ and $(z, (1+i)/2)$ and their vertical analogues also are fixed by $F$ with stabilizer containing an order four complex reflection; despite not being part of the orbifold locus for $Y$, these determine four more curves of orbifold weight four on $Z$. These are all the fixed curves for $F$, and hence determine the orbifold locus on $Z$.

From the intersection patters of all these curves, one obtains the arrangement of curves on $Z$ depicted in Figure \ref{fig:GaussA}. Note that this line arrangement on $Z$ is precisely the arrangement given by the blowup of the complete quadrangle on $\bbP^2$ at its vertices. To see this, first notice that every curve in the arrangement has self-intersection $-1$. Then, blowing down the four curves of orbifold weight eight gives $\bbP^2$, and the image of the arrangement is combinatorially the complete quadrangle. It is a classic theorem in projective geometry that every combinatorial complete quadrangle is projectively equivalent to the complete quadrangle \cite[6.13]{Coxeter}, hence the claim. Then \cite[Thm.\ 4.1(i)]{KLW} identifies the pair $(Z, \calE)$ as the Deligne--Mostow lattice with weights $(4,3,3,3,3)/8$, which is arithmetic \cite[p.\ 86]{DeligneMostow}. This completes the proof of the theorem.
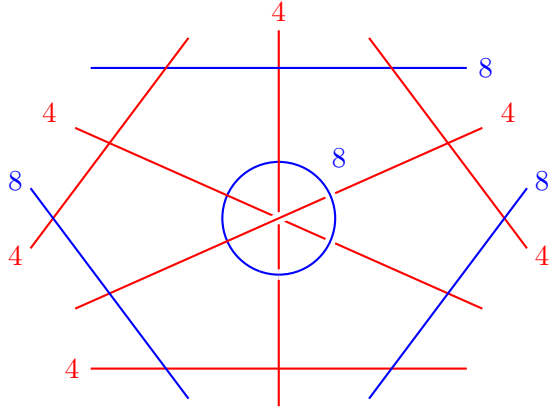
\begin{figure}[h]
\centering
\begin{tikzpicture}
\draw[thick, blue] (0,0) circle (0.75cm) node [xshift = 0.8cm, yshift=0.8cm] {$8$};
\draw[line width=4pt, white, shorten >= -0.5cm, shorten <= -0.5cm] (0,0) -- (0,-2);
\draw[thick, red, shorten >= -0.5cm, shorten <= -0.5cm] (0,2) node[yshift=0.75cm] {$4$} -- (0,-2);
\draw[line width=4pt, white, shorten >= -0.5cm, shorten <= -0.5cm] (2.25, -1) -- (0, 0);
\draw[thick, red, shorten >= -0.5cm, shorten <= -0.5cm] (2.25, -1) -- (-2.25, 1) node [xshift=-0.8cm, yshift=0.4cm] {$4$};
\draw[line width=4pt, white, shorten >= -0.5cm, shorten <= -0.5cm] (0, 0) -- (2.25, 1);
\draw[thick, red, shorten >= -0.5cm, shorten <= -0.5cm] (-2.25, -1) -- (2.25, 1) node [xshift=0.8cm, yshift=0.4cm] {$4$};
\draw[line width=4pt, white] (0,-0.75) arc (-90:45:0.75);
\draw[thick, blue] (0,-0.75) arc (-90:45:0.75);
\draw[line width=4pt, white] (0,-0.75) arc (270:250:0.75);
\draw[thick, blue] (0,-0.75) arc (270:250:0.75);
\draw[thick, blue, shorten >= -0.5cm, shorten <= -0.5cm] (-2,2) -- (2,2) node [xshift=0.75cm] {$8$};
\draw[thick, red, shorten >= -0.5cm, shorten <= -0.5cm] (-2,-2) node [xshift=-0.75cm] {$4$} -- (2,-2);
\draw[thick, red, shorten >= -0.5cm, shorten <= -0.5cm] (-1.5,2) -- (-3,0) node [xshift=-0.5cm, yshift=-0.5cm] {$4$};
\draw[thick, blue, shorten >= -0.5cm, shorten <= -0.5cm] (-1.5,-2) -- (-3,0) node [xshift = -0.5cm, yshift=0.5cm] {$8$};
\draw[thick, red, shorten >= -0.5cm, shorten <= -0.5cm] (1.5,2) -- (3,0)  node [xshift=0.5cm, yshift=-0.5cm] {$4$};
\draw[thick, blue, shorten >= -0.5cm, shorten <= -0.5cm] (1.5,-2) -- (3,0) node [xshift = 0.5cm, yshift=0.5cm] {$8$};
\end{tikzpicture}
\caption{The orbifold locus on $Z$ with weights indicated}\label{fig:GaussA}
\end{figure}
\end{pf}

\begin{rem}\label{rem:NA2}
Applying the same argument to the example mentioned in Remark \ref{rem:NA1} gives the Deligne--Mostow orbifold with weights $(6,5,5,4,4)/12$, which is \emph{nonarithmetic} \cite[p.\ 86]{DeligneMostow}.
\end{rem}

\begin{rem}\label{rem:HolzA}
One can also use the techniques from this section to prove that Holzapfel's noncompact example in \cite{HolzapfelAbelian} is commensurable with the Picard modular group over the Gaussian integers, and hence is arithmetic. Holzapfel claimed in \cite{HolzapfelAbelian} that this result would appear in a later paper, but that paper does not seem to be available.
\end{rem}

The only remaining examples from \S \ref{sec:Kodaira} that are not immediately arithmetic by virtue of having an arithmetic orbifold cover by construction are the $\bbZ / 3$ and $(\bbZ / 3)^2$ bielliptic examples contained in the proof of Theorem \ref{thm:Kodaira}. By Theorem \ref{thm:Ex2A}, it suffices to find an explicit orbifold commonly covered by the example from \S \ref{sec:Ex2} and the two bielliptic examples. Notation from the construction of each is used significantly in what follows.

It ends up being natural to start with the $(\bbZ / 3)^2$ example, since it has the least symmetry. Recall that this is constructed by a blowup $Y_9$ of the quotient $Z_9$ of $T^\prime \times T$ by the group $F_9$ generated by the automorphisms
\begin{align*}
\psi_a(z, w) &= (z + 1, \rho w) & \psi_b(z, w) =\! \left(z + \frac{\rho - 1}{3}, w + \frac{\rho - 1}{3}\right)
\end{align*}
with $\rho = \zeta^2$. Since $\tau = (\rho - 1)/3 \in T$ is one of the three fixed points of $\rho$ acting on $T$, the map $r(z, w) = (z, \rho w)$ commutes with $\psi_a$ and $\psi_b$. The first lemma gives a birational model for the eventual common orbifold quotient.

\begin{lem}\label{lem:Z9quo}
The automorphism $r$ of $T^\prime \times T$ induces a nontrivial order three automorphism $r_9$ of $Z_9$ with quotient isomorphic to the quotient $B_0$ of $T \times \bbP^1$ by the order three automorphism $(z, y) \mapsto (z + \tau, \rho y)$.
\end{lem}

\begin{pf}
That $r$ induces an automorphism of $Z_9$ follows from the fact that $r$ commutes with the generators $\psi_a, \psi_b$ of $F_9$. Since $r$ has fixed points, and thus is not in the deck group $F_9$, the induced automorphism of $Z_9$ is nontrivial.

The goal is then to understand the quotient of $T^\prime \times T$ by the $(\bbZ / 3)^3$ group of automorphisms generated by $\psi_a$, $\psi_b$, and $r$, which can be alternately generated by $r$ along with $r^{-1} \circ \psi_a$ and $\psi_b$, which both act by translations. The quotient of $T^\prime \times T$ by $r^{-1} \circ \psi_a$ is $T \times T$ by definition and the transformation of $T \times T$ induced by $\psi_b$ is diagonal translation by $\tau$. Note that the quotient of $T$ is by the action of translation by $\tau$ is $(\rho - 1)^{-1} \bbZ[\zeta] \bs \bbC$, which is homothetic and thus isomorphic to $T$. The action on $T \times T$ induced by the action of $r$ is the action by $\rho$ on the second coordinate, so the quotient is $T \times \bbP^1$. Thus there is a commutative diagram
\[
\begin{tikzcd}[column sep = small, row sep = small]
T^\prime \times T \arrow[dd] \arrow[r] & Z_9 & & & & \\
& & & & & \\
T \times T \arrow[dd] \arrow[rrrrr] & & & & & T \times \bbP^1 \arrow[dd] \\
& & & & & \\
A = (T \times T) / \langle (\tau, \tau) \rangle \arrow[dd] \arrow[rrrrr] & & & & & B_0 \arrow[dd] \arrow[from=lllluuuu, crossing over] \\
& & & & & \\
T \times T \arrow[rrrrr] & & & & & T \times \bbP^1
\end{tikzcd}
\]
where the bottom $T \times T$ is the quotient of $T \times T$ by the $(\bbZ / 3)^2$ generated by the action of $\tau$ on each individual coordinate, $T^\prime \times T \to Z_9$ is the covering projection, and every other right-moving arrow is the action induced by $r$.

All that remains is to understand the automorphism of $T \times \bbP^1$ induced by translation by $(\tau, \tau)$ on $T \times T$. Since $r$ only changes the second coordinate, the action on the first factor is evidently by $\tau$. A point $y \in \bbP^1$ with preimage $z \in T$ is fixed under the induced action if and only if $\rho z = z + k \tau$ for some $k \in \{0,1,2\}$, i.e., $(\rho - 1) z = k \tau$. This implies the induced automorphism of $\bbP^1$ is nontrivial of order three, and hence is multiplication by $\rho$ in the appropriate coordinates. This proves the lemma.
\end{pf}

The next goal is to see that $r_9$ induces an automorphism of the ball quotient associated with the pair $(Y_9, \calD_9)$.

\begin{lem}\label{lem:Worb9}
The action of $r_9$ on $Z_9$ induces an action on the blowup $Y_9$. Let $B$ denote the normal surface $Y_9 / \langle\wh{r}_9\rangle$ completing the commutative diagram
\[
\begin{tikzcd}[column sep = small, row sep = small]
Z_9 \arrow[d] & Y_9 \arrow[l] \arrow[d] \\
B_0 & B \arrow[l]
\end{tikzcd}
\]
with birational horizontal arrows. Then there is a $\bbQ$-divisor $\calE$ on $B$ so that $(B, \calE)$ is a ball quotient pair with $(Y_9, \calD_9)$ as an orbifold cover.
\end{lem}

\begin{pf}
The first thing to prove is that $r_9$ induces an automorphism of the blowup $Y_9$. Let $D$ denote the image in $T^\prime \times T$ of the diagonal of $\bbC^2$, which projects to the curve $C_0^\prime$ in $Z_9$. Notice that $\psi_b$ is diagonal translation and thus fixes $D$ and that $r(D) = \psi_a(D)$ since $\rho z = \rho(z + 1)$ on $T$. This implies that $r$ preserves the $F_9$-orbit of $D$ and thus $r_9$ stabilizes $C_0^\prime$. To fix $C_0^\prime$, $r_9$ must fix its unique triple point, and thus $r_9$ acts on the blowup $Y_9$ by a transformation $\wh{r}_9$ that fixes the proper transform $C_0$ of $C_0^\prime$.

Note that the fixed point $(0,0) \in D \subset T^\prime \times T$ of $r$ maps to the singular point of $C_0^\prime$. The derivative of $r$ at $(0,0)$ is not a scalar multiple of the identity, which implies that the action of $\wh{r}_9$ on the exceptional divisor $E$ of $Y_9$ is nontrivial, i.e., $\wh{r}_9$ does not act by a complex reflection through the curve. On the other hand, $r$ does act as a complex reflection through the lift to $T^\prime \times T$ of the horizontal curve $C_1^\prime$ through the singular point of $C_0^\prime$, namely the curve with coordinates $(z, 0)$, hence $\wh{r}_9$ also acts as a complex reflection through the proper transform $C_1$ of $C_1^\prime$ to $Y_9$. Thus $\wh{r}_9$ fixes each curve in the orbifold locus for the pair $(Y_9, \calD_9)$ and the quotient $B$ of $Y_9$ by $\langle \wh{r}_9 \rangle$ is a ball quotient orbifold.
\end{pf}

The next lemma completely describes the $\bbQ$-divisor for the pair $(B, \calE)$.

\begin{lem}\label{lem:Z9E}
The curves in the orbifold locus for the pair $(B, \calE)$ are:
\begin{itemize}

\item[$\star$] the image $\wh{E}$ of the exceptional curve $E$ of $Y_9$, with orbifold weight $3$

\item[$\star$] the image $\wh{C}_0$ of the curve $C_0$ on $Y_9$, with orbifold weight $3$

\item[$\star$] the image $\wh{C}_1$ of the horizontal curve $C_1$ on $Y_9$, with orbifold weight $9$

\end{itemize}
\end{lem}

\begin{pf}
From the discussion before the statement of the lemma, it suffices to show that $C_1$ is the only curve on $Y_9$ on which $\wh{r}_9$ has a fixed point. Indeed, the given curves were shown to be contained in the support of $\calE$ with the given weight, so it suffices to show that $\wh{r}_9$ does not have either an isolated fixed point or another fixed curve on which it acts by a complex reflection. Suppose that $(z,w) \in T^\prime \times T$ projects to a fixed point of $\wh{r}_9$. Then there is an element $g$ of the covering group for $Z_9$ so that $g(z,w) = (z, \rho w)$. In other words, there are $k, \ell \in \{0,1,2\}$ so that
\[
(z, \rho w) =\! \left(z + k + \ell \frac{\rho - 1}{3}\,,\, \rho^k w + \ell \tau\right)
\]
on $T^\prime \times T$. Then evidently $k = \ell = 0$, and thus $\rho w = w$. In other words, the only such curves are $(z,0)$, $(z, \tau)$, and $(z, 2 \tau)$, which all map to the horizontal curve $C_1^\prime$ with proper transform $C_1$. This proves the lemma.
\end{pf}

\begin{prop}\label{prop:Z39comm}
The ball quotient pair $(B, \calE)$ described by Lemma \ref{lem:Z9E} is also a quotient of the $\bbZ/3$ bielliptic example from \S \ref{sec:Kodaira}, hence the $\bbZ / 3$ and $(\bbZ / 3)^2$ bielliptic ball quotient orbifolds are commensurable.
\end{prop}

\begin{pf}
The $\bbZ / 3$ bielliptic example is constructed using the automorphism
\[
\phi_3(z, w) = (z + \tau, \rho w)
\]
of $T \times T$, which commutes with the automorphism $t(z, w) = (z + \tau, w + \tau)$ of order three. Moreover, $t$ permutes the three points $(\rho^j/3, \rho^j/3)$, $j \in \{0,1,2\}$, that project to the triple point on the curve $C_0^\prime$ on the bielliptic surface $Z_3$. It follows that $t$ induces an automorphism $\wh{r}_3$ of the blowup $Y_3$ of $Z_3$.

To compute the action of $\wh{r}_3$ on the exceptional divisor of $Y_3$, one can consider the derivative of $\phi_3^{-1} \circ t$ fixing the point $(1/3, 1/3)$ on $T \times T$. This is the automorphism
\[
(z,w) \mapsto (z, \rho^2 w + \tau),
\]
which has nonscalar derivative at $(1/3, 1/3)$. This implies that the action of $\wh{r}_3$ stabilizing the exceptional curve is not by a complex reflection.

Any other fixed point of $\wh{r}_3$ is determined by a fixed point of $\phi_3^j \circ t$ for some $j \in \{0,1,2\}$. The equation
\[
(z,w) =\! \left(z + \tau + j \tau, \rho^j w + \tau\right),
\]
implies that $z = z + (1 + j) \tau$, so $j = 2$. Then $\rho^2 w + \tau = w$, which has only the three solutions $w = \rho^k / 3$ for $k \in \{0,1,2\}$. In other words, the only fixed points of $\wh{r}_3$ on $Z_3$ are the horizontal curve $C_1$ through the exceptional divisor, and $\wh{r}_3$ acts on this curve as a complex reflection.

As in many previous cases, to show that $\wh{r}_3$ acts on the ball quotient pair $(Y_3, \calD_3)$ it only suffices to show that $\wh{r}_3$ stabilizes the curve $C_0$ on $Y_3$ determined by the diagonal of $T \times T$. However, $t$ stabilizes the diagonal of $T \times T$ by definition, so $\wh{r}_3$ indeed stabilizes $C_1$ and thus fixes all the curves in the support of $\calD_3$.

Therefore the proof of the proposition will be complete once it is shown that the quotient ball quotient pair is the pair $(B, \calE)$ from Lemma \ref{lem:Z9E} and that $\calD_3$ maps to $\calE$ under the quotient map, with $C_1$ mapping to the curve of orbifold weight $9$. One of the equivalent definitions of $B_0$ given in the proof of Lemma \ref{lem:Z9quo} is by taking the quotient of $T \times T$ first by the diagonal action of $(\tau, \tau)$ (i.e., the transformation $t$) followed by the action of $\rho$ on the second coordinate, which is precisely the induced action of $\phi_3$ on the quotient $(T \times T) / \langle (\tau, \tau) \rangle$. Moreover, this description implies that $\calD_3$ indeed maps to $\calE$ with the appropriate orbifold weights. Thus $(Y_3, \calD_3)$ and $(Y_9, \calD_9)$ are commensurable.
\end{pf}

\begin{thm}\label{thm:Z39A}
The $\bbZ/3$ and $(\bbZ / 3)^2$ bielliptic ball quotients constructed in the proof of Theorem \ref{thm:Kodaira} are arithmetic.
\end{thm}

\begin{pf}
Since Proposition \ref{prop:Z39comm} shows that both bielliptic examples are commensurable with the ball quotient pair $(B, \calE)$ from Lemma \ref{lem:Z9E}, it suffices to prove arithmeticity of that example. The proof will show that the example from \S\ref{sec:Ex2} admits a $\bbZ / 3$ action with quotient $(B, \calE)$. Accordingly, this proof freely draws from the notation established in \S \ref{sec:Ex2} and \S \ref{sec:Arith2}.

Consider the automorphism $b = \rho^2 \beta$ of $T \times T$, which has matrix
\[
\begin{pmatrix} \rho^2 & 0 \\ \rho^2 & 1 \end{pmatrix}
\]
and eigenvectors $(1-\rho, 1)$ and $(0, 1)$ with respective eigenvalues $\rho^2$ and $1$. It is already known from the proof of Proposition \ref{prop:FoundDM} that $b$ acts as an automorphism of the ball quotient pair $(Y, \calD)$. It acts nontrivially (i.e., not as a complex reflection) on the exceptional curve, cyclically permutes the three curves $(z, 0)$, $(z, z)$, and $(z, \zeta z)$, and acts on the curve $(0, z)$ as a complex reflection. The only other curve fixed by $b$ is the proper transform of the eigenline $((1 - \rho) z, z)$, where $b$ acts as complex multiplication by $\rho^2$, in particular not by a complex reflection. It follows that the orbifold locus for the quotient of $Y$ by $\langle b \rangle$ has the same combinatorics and orbifold weights as $\calE$. It remains to prove that $Y / \langle b \rangle$ is isomorphic to $B$ and that the orbifold locus $\calD$ appropriately maps to $\calE$ under the projection.

Decompose $T \times T$ in terms of the two eigenlines $((1-\rho) x, x)$ and $(0, y)$ of $b$. This decomposition is nonunique. Specifically, every point $(z, w)$ can be written in coordinates as
\[
(z, w) = ((1 - \rho) x, x) + (0, y)
\]
where $y$ is uniquely defined and $x$ is only defined up to translation by $\langle \tau \rangle$. In these coordinates, the action of $b$ is by $(x, y) \mapsto (\rho^2 x, y)$. More concretely, consider the self-covering of $T \times T$ by itself induced by $\tau$ acting on the first coordinate. Then $b$ is the automorphism of the quotient induced by $\rho^2$ acting on the first coordinate of the cover. The quotient of the cover $T \times T$ by the lift of $b$ is then $\bbP^1 \times T$, and the further quotient by the automorphism induced by the action of $\tau$ is isomorphic to $B_0$ by construction of $B_0$. Considering the curve stabilizer calculations for the $b$-action, it follows that the quotient orbifold is indeed the pair $(B, \calE)$. This proves the theorem.
\end{pf}

This finally completes all the pieces required to prove Theorem \ref{thm:Arithmetic}.

\begin{pf}[Proof of Theorem \ref{thm:Arithmetic}]
By construction, all ball quotient orbifolds used to prove Theorems \ref{thm:Main} and \ref{thm:Kodaira} are explicitly commensurable with the orbifolds shown to be arithmetic in Theorems \ref{thm:Ex1A}, \ref{thm:Ex2A}, \ref{thm:GaussExA}, and \ref{thm:Z39A}.
\end{pf}

\bibliography{Wiman}

\end{document}